\documentclass[a4paper,12pt,twoside]{article}

\usepackage[british]{babel}
\usepackage{url}

\usepackage{amsmath,amssymb,amsfonts} \usepackage{epsfig, color,
  subfig} \usepackage{enumerate,booktabs,makecell}
\usepackage{bm}

\newtheorem{theorem}{Theorem}[section]
\newtheorem{lemma}[theorem]{Lemma}

\def\vet#1{{\bm #1}}
\def\build#1_#2^#3{\mathrel{\mathop{\kern 0pt#1}\limits_{#2}^{#3}}}

\def\reali{\mathbb{R}}
\def\complessi{\mathbb{C}}

\def\interi{\mathbb{Z}}
\def\toro{\mathbb{T}}

\def\Ascr{\mathcal{A}}
\def\Bscr{\mathcal{B}}

\def\Escr{\mathcal{E}}
\def\Fscr{\mathcal{F}}
\def\Gscr{\mathcal{G}}
\def\Lbf{\mathbf{L}}
\def\Hscr{\mathcal{H}}
\def\Jscr{\mathcal{J}}
\def\Kscr{\mathcal{K}}
\def\Lscr{\mathcal{L}}
\def\Nscr{\mathcal{N}}
\def\Oscr{\mathcal{O}}
\def\Pscr{\mathcal{P}}
\def\Rscr{\mathcal{R}}
\def\Sscr{\mathcal{S}}

\def\epsilon{\varepsilon}
\def\rho{\varrho}
\def\phi{\varphi}

\def\Pgot{{\mathfrak P}}

\def\lie#1{\Lscr_{#1}}
\def\Lie#1{\Lscr_{#1}}
\def\imunit{{\bf i}}
\def\poisson#1#2{\lbrace #1,#2 \rbrace}
\def\fastpoisson#1#2{\left\{ #1,#2 \right\}_{\vet{L},\vet{\lambda}}}
\def\secpoisson#1#2{\left\{ #1,#2 \right\}_{\vet{\xi},\vet{\eta}}}

\def\Pgot{{\mathfrak P}}

\def\wideitem#1{\par\hangindent\itemindent
   \noindent\hbox to\parindent{\hfil{#1}\enspace}\ignorespaces}

\title{\bf Librational KAM tori in the secular dynamics of the
  $\upsilon$~Andromed{\ae} planetary system\thanks{\noindent{\it Key
      words and phrases:} exoplanets, Celestial Mechanics.  
}
}

\author{
{\bf CHIARA CARACCIOLO}\\
{\small Dipartimento di Matematica, Universit\`a degli Studi di Milano,}\\
{\small via Saldini 50, 20133\ ---\ Milano, Italy.}\\
{\bf UGO LOCATELLI}\\
{\small Dipartimento di Matematica, 
Universit\`a degli Studi di Roma ``Tor Vergata'',}\\
{\small via della Ricerca Scientifica 1, 00133\ ---\ Roma, Italy.}\\
{\bf MARCO SANSOTTERA}\\
{\small Dipartimento di Matematica, Universit\`a degli Studi di Milano,}\\
{\small via Saldini 50, 20133\ ---\ Milano, Italy.}\\
{\bf MARA VOLPI}\\
{\small Dipartimento di Matematica, 
Universit\`a degli Studi di Roma ``Tor Vergata'',}\\
{\small via della Ricerca Scientifica 1, 00133\ ---\ Roma, Italy.}\\
{\small e-mails:
  {\tt chiara.caracciolo@unimi.it, locatell@mat.uniroma2.it,}}\\
{\small {\tt marco.sansottera@unimi.it, volpi@mat.uniroma2.it}}
}

\begin{document}
\maketitle


\markboth{C. Caracciolo, U. Locatelli, M. Sansottera, M. Volpi}{Librational KAM tori in the secular dynamics of $\upsilon$~Andromed{\ae} planetary system}

\begin{abstract}
  \noindent
   We study the planetary system of
   $\upsilon$~Andromed{\ae}, considering the three-body
   problem formed by the central star and the two largest planets,
   $\upsilon$~And~\emph{c} and
   $\upsilon$~And~\emph{d}.  We adopt a secular,
   three-dimensional model and initial conditions within the range of
   the observed values. The numerical integrations highlight that the
   system is orbiting around a one-dimensional elliptic torus (i.e., a
   periodic orbit that is linearly stable). This invariant object is
   used as a seed for an algorithm based on a sequence of canonical
   transformations. The algorithm determines the normal form related
   to a KAM torus, whose shape is in excellent agreement with the
   orbits of the secular model. We rigorously prove that the algorithm
   constructing the final KAM invariant torus is convergent, by
   adopting a suitable technique based on a computer-assisted proof.
\end{abstract}


\section{Introduction}
\label{sec:intro}
The discovery of the first multiple-planet extrasolar system dates
back to the very end of the last century
(see~\cite{Butler-et-al-1999}). Through the radial velocity method,
three planets were found orbiting around one of the two stars of the
binary system $\upsilon$~Andromed{\ae}. Since then thousands of
exoplanets have been observed, and among the plethora of new
extrasolar planetary systems around $780$~of them host two or more
planets\footnote{Statistical data on the confirmed detections of
  exoplanets can be easily produced by accessing to the website {\tt
    http://exoplanet.eu/}}.  Of all the different detection techniques
available at present day, the radial velocity method better highlights
the skeleton of a multi-planetary system, as it is particularly
sensitive to more massive bodies. The dynamical characteristics of
radial-velocity-detected systems can strongly differ from those of the
solar system. For instance, the orbits can be far more eccentric and
the planets can present masses several times larger than the one of
Jupiter. A proper understanding of the architecture of these systems
is complicated even more by the fact that some of their orbital
parameters are unknown. Concerning the radial velocity method, this
applies in particular to the orbital inclination $i_P\,$ of an
exoplanet with respect to the line of sight (see,
e.g.,~\cite{Beau-FerM-Mich-2012}). This implies that only minimal
masses for the planets can be inferred (i.e., $m_P\sin i_P$), whereas
their actual values remain unknown.

In several cases, exoplanets can be observed by only one detection
technique. Therefore, it is a rare occurrence that two different
detection methods can be used on the same system.  The planetary
system of $\upsilon$~Andromed{\ae} is then particularly interesting
because it was thoroughly observed by means of both radial velocity
and astrometry. This makes $\upsilon$~Andromed{\ae} one of the most
well-known and studied exoplanetary systems. Thanks to the joint
observations, it was possible to determine ranges of values for the
orbital inclinations of the two planets with the largest minimal
masses, $\upsilon$~And~\emph{c} and $\upsilon$~And~\emph{d}
(see~\cite{McArt-et-al-2010}). Therefore, ranges of possible values
for both planetary masses have been inferred as well. These new
constraints modify substantially the study of the stability of the
system. In fact, since the orbital plane of $\upsilon$~And~\emph{c} is
highly inclined with respect to the line of sight, the actual value of
the mass could be larger than $5$~times the minimal
one\footnote{According to the data reported in~\cite{McArt-et-al-2010}
  the uncertainty on that measure is still relevant, being its
  half-width about 30\% of the mid value.}.  The same effect has a
weaker impact on the determination of $\upsilon$~And~\emph{d}'s mass,
being the increasing factor about $2.5\,$.  The two-dimensional models
of the $\upsilon$~Andromed{\ae} system that were studied mainly in the
first decade of this century used the minimal values of the masses.
Therefore, the perturbation of the Keplerian orbits (that is due to
the mutual gravitation) when considering the updated values is
expected to be one order of magnitude larger. Moreover, the mutual
inclination between $\upsilon$~And~\emph{c} and
$\upsilon$~And~\emph{d} is quite remarkable, being
about~$30^\circ$. Taking into account that also their eccentricities
are about 5~times larger than those we observe in our Jovian planets,
it is natural to expect that the excursions of the secular variables
can be very relevant, with a possibly dramatic impact on the stability
of the system.  This is confirmed by the fact that, according to the
analysis described\footnote{Indeed, in~\cite{Deitrick-et-al-2015} a
  model including $\upsilon$~And~\emph{b}, $\upsilon$~And~\emph{c} and
  $\upsilon$~And~\emph{d} was considered. However, the instability for
  many initial conditions (that are compatible with the observations)
  can be easily confirmed by numerical integrations studying
  separately the dynamics of the two latter planets, see, e.g.,
  C. Caracciolo: {\it On the stability in the neighborhood of
    invariant elliptic tori}, Ph.D. thesis, Univ. of Rome ``Tor
  Vergata'' (2021). We also recall that in~\cite{McArt-et-al-2010} a
  long-period trend indicating the presence of a fourth planet in the
  system was also detected. However, in the present work we will study
  a model where $\upsilon$~And~\emph{b} and $\upsilon$~And~\emph{e}
  are completely neglected.\label{phdtesi-Chiara}}
in~\cite{Deitrick-et-al-2015}, the orbits starting from many initial
conditions within the range of observed values are evidently
unstable. Therefore, the stability of the $\upsilon$~Andromed{\ae}
planetary system in a genuine 3D-model cannot be ensured by simply
applying the AMD general criterion (see~\cite{Las-Pet-2017}, where
this method was introduced).

This paper aims to study the stability of the $\upsilon$~Andromed{\ae}
system in the refined framework of the KAM theory
(see~\cite{Kolmogorov-1954}, \cite{Moser-1962}
and~\cite{Arnold-1963}). Our approach is definitely challenging from a
computational point of view, but we will show that it is able to
construct tori that are invariant with respect to the secular model of
the planetary Hamiltonian. In order to adapt our KAM-tori-constructing
method as it has already been applied to other problems in Celestial
Mechanics (see, e.g.,~\cite{Loc-Gio-2000} and~\cite{Gab-Jor-Loc-2005})
to this particular context, a sequence of canonical transformations
must be performed as a preparatory step.  In particular, we will
stress the key role that is played in our approach by the preliminary
construction of the normal form of a so called 1D elliptic torus
(i.e., a periodic orbit that is linearly stable). The final
application of a suitable adaptation of the classical Kolmogorov
normalization algorithm works because such a periodic orbit is not far
from the initial conditions that we consider within the range of the
observations.  In fact, this distance (in the phase space) is so small
that the motion we study remains trapped in a resonance such that the
difference of the arguments of pericenter is in libration. In this
respect, our work relies also on some ideas found in the studies that
in the first decade of this century were devoted to the dynamical
orbits of multi-planet extrasolar systems. In fact, just a couple of
years after their discovery $\upsilon$~And~\emph{c} and
$\upsilon$~And~\emph{d} were conjectured to be into a librational
regime with respect to the difference of the pericenters arguments
(see~\cite{Chi-Tab-Tre-2001}).  Moreover, three years later the
relevance of such an apsidal locking was emphasized into the framework
of the secular model representing the slow dynamics of three--body
planetary systems with rather eccentric orbits
(see~\cite{Mich-Mal-2004}).

Our paper is organized as follows. In Section~\ref{sec:model}, we
introduce the secular Hamiltonian model of the $\upsilon$~Andromed{\ae}
planetary system. In Section~\ref{sec:toro} we fully describe the
algorithm constructing invariant tori (of both elliptic and KAM type)
for such a model. The obtained results are discussed in
Section~\ref{sec:results}. The conclusions are outlined in
Section~\ref{sec:conclu}. Furthermore, some technical details about
the rigorous computer-assisted proof of the existence of KAM tori are
deferred to Appendix~\ref{sec:technicalities}.

\section{Settings for the definition of the Hamiltonian model}
\label{sec:model}

\subsection{Planetary three-body problem}
\label{sbs:3body-model}
We consider the planetary three--body problem, i.e., we focus on an
Hamiltonian model including two planets orbiting around a star under
the mutual effect of the gravitational forces.
Let $P_0$ be the star and $P_1$ and $P_2$ be the two planets having
mass $m_j$ (with $j=0,1,2$).  The problem, as we are in the spatial
case, has $9$ degrees of freedom, which can be reduced to $6$ by using
the conservation of the linear momentum. If $(\vet r_1, \vet r_2,
\tilde {\vet r}_1, \tilde {\vet r}_2)$ are the canonical variables
with respect to an heliocentric reference frame, that are $\vet r_1 =
\overrightarrow{P_0 P_1}$, $\vet r_2 = \overrightarrow{P_0 P_2}$ and
$\tilde{\vet r}_1$, $\tilde{\vet r}_2$ the corresponding momenta, then
the Hamiltonian of the problem takes the form
\begin{equation}
\label{frm:ham3cpiniz}
H = \frac 1 2 \sum_{i=1}^2 \| \tilde {\vet r}_i\|^2
\left[\frac{1}{m_i}+ \frac{1}{m_0}\right] - G \sum_{i=1}^2 \frac{m_0
  m_i}{r_i} +\sum_{0<i<j} \frac{\tilde{\vet r}_i\cdot \tilde{\vet
    r}_j}{m_0} - G \sum_{0< i<j} \frac{m_i m_j}{\Delta_{ij}}\ ,
\end{equation}
where $G$ is the gravitational constant, $\Delta_{ij} = \|\vet r_i -
\vet r_j\|$ and we have distinguished between the kinetic energy and
the potential one for what concerns both the contributions given by
the Keplerian part and the one due to the interactions between the
planets. The degrees of freedom can be further reduced by 2 using the
conservation of the total angular momentum $\vet C$. This allows us to
write the Hamiltonian in Poincar\'e variables, that are
\begin{align}
\label{frm:poin-var}
\Lambda_j &= \frac{m_0 m_j}{m_0 + m_j} \sqrt{G(m_0+m_j)a_j}\ , &\lambda_j
=& M_j + \omega_j\ ,\\ \xi_j & = \sqrt{2 \Lambda_j}
\sqrt{1-\sqrt{1-e_j^2}}\cos{(\omega_j)}\ , &\eta_j =& - \sqrt{2
  \Lambda_j} \sqrt{1-\sqrt{1-e_j^2}}\sin{(\omega_j)}\ ,
\end{align}
where $a_j$, $e_j$, $M_j$, and $\omega_j$ are the semi-major axis, the
eccentricity, the mean anomaly and the argument of the pericenter of
the planet $j$, respectively.  The reduction of the total angular
momentum makes implicit the dependence on the inclinations $i_j$ and
on the longitudes of the nodes $\Omega_j$. For the sake of clarity,
let us report here the value of the mutual inclination, that in the
Laplace reference frame is the sum of the two inclinations, i.e.,
\begin{equation}
i_1 +i_2 = \arccos\left(\frac{C^2-\Lambda_1^2(1-e_1^2)-\Lambda_2^2(1-e_2^2) }{2\Lambda_1 \Lambda_2 \sqrt{1-e_1^2}\sqrt{1-e_2^2}}\right)\ ,
\end{equation}
being $C$ the module of the total angular momentum 
\begin{equation}
\label{frm:ang-mom}
 C = \sum_{k=1}^2 \Lambda_k \sqrt{1-e_k^2}\cos i_k\ .
\end{equation}
Moreover, we introduce a translation $L_j = \Lambda_j - \Lambda_j^*$,
where $\Lambda_j^*$ is defined in order to obtain that in the
Keplerian approximation of the motion the values of the semi-major
axes are in agreement with the observations. In view of the Poisson
theorem (see~\cite{Dur-1978}), it is known that the semi-major axes
have not any secular contribution up to order two in the masses. This 
is the reason why the expansions are usually made around the average
values of the semi-major axes or their initial values. For the sake of
simplicity, we will adopt this latter option.

The Hamiltonian is expanded with respect to these Poincar\'e
variables\footnote{The computation of this Hamiltonian starting from
  the one in~\eqref{frm:ham3cpiniz} is not straightforward. The main
  difficulty concerns the expansion of the inverse of the
  distance between the planets $\Delta_{ij}$. For a detailed
  discussion of the method used here for doing such a calculation we
  defer again to~\cite{Las-1989}.} and the parameter $D_2$, that
measures the difference between the total angular momentum of the
system and the one of a similar system with circular and coplanar
orbits; i.e., it is defined as
\begin{equation}
\label{frm:D2}
D_2 = \frac{(\Lambda_1^*+\Lambda_2^*)^2 -C^2}{\Lambda_1^* \Lambda_2^*};
\end{equation}
therefore, it is of the same order as $e_1^2+i_1^2+e_2^2+i_2^2\,$.
Thus, we can write the Hamiltonian of the three--body problem as
\begin{equation}
 \label{frm:3BP}
  H= \Kscr(\vet L) + \Pscr(\vet L, \vet \lambda, \vet \xi, \vet \eta;D_2) = \sum_{j_1=1}^{\infty}h^{({\rm Kep})}_{j_1,0}(\Lbf) +
  \mu\sum_{s=0}^{\infty}\sum_{j_1=0}^{\infty}\sum_{j_2=0}^{\infty}\,D_2^s\,
  h^{(\mathcal{P})}_{s;j_1,j_2}(\Lbf,\boldsymbol{\lambda},\boldsymbol{\xi},\boldsymbol{\eta})
\end{equation}
where $\mu = \max\{m_1/m_0,m_2/m_0\}$. Moreover,
\begin{itemize}
\item $ \Kscr(\vet L)= \sum_{j_1=1}^{\infty}h^{({\rm Kep})}_{j_1,0}(\Lbf)$
  is the Keplerian part and $h^{({\rm Kep})}_{j_1,0}$ is a
  homogeneous polynomial of degree $j_1$ in $\Lbf$; in
  particular, $h^{({\rm Kep})}_{1,0} = \vet{n}^*\cdot\Lbf$, where the
  components of the angular velocity vector $\vet{n}^*$ are defined by
  the third Kepler law;
\item$h^{(\Pscr)}_{s;j_1,j_2}$ is a homogeneous polynomial of
  degree $j_1$ in $\Lbf$, degree $j_2$ in
  $(\boldsymbol{\xi},\boldsymbol{\eta})$,
  with coefficients that are trigonometric polynomials in
  $\boldsymbol{\lambda}$, being the harmonics $\vet k \cdot \vet
  \lambda$ such that $|\vet k| = \sum_{j=1}^2 |k_j|\le s$.
\end{itemize}
Clearly, in the applications we deal with finite expansions; the
truncation parameters will be discussed when appropriate.

\subsection{Secular model at order two in the masses}
\label{sec:secular-model}
The expression of the Hamiltonian of the three-body problem
in~\eqref{frm:3BP} highlights the distinction between the so called
{\it fast variables} $(\vet L, \vet \lambda)$ and the {\it secular
  variables} $(\vet \xi, \vet \eta)$. Indeed, if we consider the
corresponding Hamilton equations, we have that $\dot {\vet \lambda} =
\Oscr(1)$. This means that the motion of the planet along the orbit,
that is in first approximation a Keplerian ellipse, has a different
timescale with respect to the secular variables, whose variation is
due to the interaction between the planets and, therefore, is of
$\Oscr(\mu)$. Since we are interested in the study of the long-time
stability of the system, a common procedure consists on considering
just the evolution of the secular variables, by averaging the
Hamiltonian with respect to the fast angles $ \vet \lambda$. With a
simple average of $\Pscr$ we would obtain a secular approximation with
terms of order $\mu$, namely at order $1$ in the masses. In this work,
we consider terms up to order $2$ in the masses, averaging with a
near-the-identity canonical change of coordinates inspired by the
algorithm for the construction of the Kolmogorov normal form. Indeed,
we focus on the torus corresponding to $\boldsymbol{L}=0$. The
procedure is quite standard in Celestial Mechanics but, for the sake
of completeness, we sketch here the main steps (see for more
details~\cite{Vol-Loc-San-2018}, that is in turn an adaptation of the
approach developed in~\cite{Loc-Gio-2000}).

The first transformation of coordinates that we define aims at
removing the perturbative terms that depend on the angles $\vet
\lambda$ but do not depend on the actions $\vet L$, being
$\dot{L_j}=\partial H/\partial\lambda_j$ for $j=1,2\,$.  This is done
by using the term linear in the actions, i.e., $\vet n^* \cdot \vet
L$, to define a generating function $\chi^{(\Oscr 2)}_1(\vet \lambda)$
as the solution of the following homological equation:
\begin{equation}
  \label{eq:homeq}
\left\{\chi^{(\Oscr 2)}_1,\ \vet n^* \cdot \vet L\right\} +\mu \sum_{{s=0\,,\>j_2=0}\atop{2s+j_2\le
      N_S}}\left\lceil D_2^s\, h_{s;0,j_2}^{(\Pscr)}\right\rceil_%
      {\boldsymbol{\lambda}:K_F}= \mu
      \sum_{{s=0\,,\>j_2=0}\atop{2s+j_2\le N_S}}D_2^s\, \Big\langle
      h_{s;0,j_2}^{(\Pscr)}\Big\rangle_{\boldsymbol{\lambda}}\ ,
\end{equation}
being $\langle \cdot \rangle_{\vet \lambda}$ the average with respect
to the angles $\vet \lambda$, while with the notation $\lceil \cdot
\rceil_{K_F}$ we mean that the expansions are truncated at the
trigonometrical degree $K_F$ in the angles $\vet \lambda$.  Let us add a
few comments about the truncations parameters $K_F$ and $N_S$. The
value of $K_F$ is defined so as to take into account the main mean-motion quasi-resonances of the system considered.  For example, if the
system is close to the resonance $k_1^*:k_2^*$, then $K_F$ is defined
as $K_F \ge |k_1^*|+|k_2^*|$. In the same spirit, the value $N_S$ of
the truncation of the expansions in eccentricity and inclination is
set in order to consider the quasi-resonance. Let us assume that the
quasi-resonant angular terms are of type $(k_1^*\lambda_1 - k_2^*
\lambda_2)$, then in view of the D'Alembert rules it is convenient up
to consider expansions up to order in eccentricity and inclination
$N_S\ge 2(|k_1^*|-|k_2^*|)$. A more accurate discussion about the
choice of these parameters is deferred to the following subsection.

Now we have to apply the transformation of coordinates defined by the
application of the Lie series operator $\exp(\Lie{\chi_1^{(\Oscr 2)}})\ \cdot = \sum_{j=0}^\infty (1/j!)\Lie{\chi_1^{(\Oscr 2)}}^j\,\cdot$ to the
Hamiltonian. Recalling that in our secular model we will not consider
terms depending on $\vet L$ or of order greater than $\mu^2$, the only
terms we need to compute are included in the following expansion:
\begin{equation}
  \label{eq:htilde}
  \begin{aligned}
    \widetilde{\kern-2pt H} &= H +
    \frac{1}{2}\left\{\chi^{(\Oscr 2)}_1, \ \Lscr_{\chi^{(\Oscr  2)}_1}h^{({\rm Kep})}_{2,0}\right\}_{\vet L, \vet \lambda}\\ &+\left\{\chi^{(\Oscr
        2)}_1,\ \mu \sum_{{s=0\,,\>j_2=0}\atop{2s+j_2\le N_S}}
      D_2^s\,{h}_{s;1,j_2}^{(\Pscr)}\right\}_{\vet L, \vet \lambda}+
    \frac{1}{2}\left\{\chi^{(\Oscr 2)}_1,\ \mu
      \sum_{{s=0\,,\>j_2=0}\atop{2s+j_2\le
          N_S}}D_2^s\,{h}_{s;0,j_2}^{(\Pscr)}\right\}_{\vet \xi,\vet \eta}
  \end{aligned}
\end{equation}
where $\fastpoisson{\cdot}{\cdot}$ and $\secpoisson{\cdot}{\cdot}$ are
the terms of the Poisson bracket involving only the derivatives with
respect to the pairs of conjugate variables $(\vet L,\vet \lambda)$ and $(\vet \xi,\vet \eta)$,
respectively. Then, according to~\cite{Loc-Gio-2000}, we have that 
\begin{equation*}
  \langle H^{(\Oscr 2)} \rangle_{\vet{\lambda}} \Big|_{\Lbf= \vet{0}}
  = \langle\, \widetilde{\kern-2pt H}\, \rangle_{\vet{\lambda}}
  \Big|_{\Lbf= \vet{0}} + {\cal O} (\mu^3)\,,
\end{equation*}
being $H^{(\Oscr 2)} = \exp(\Lie{\chi_1^{(\Oscr2)}})H$.  Let us remark
that for the definition of this model it is not necessary to compute
the effects induced by the second generating function
$\chi_2^{(\Oscr_2)}(\vet L,\vet \lambda)$ for removing terms linear in
$\vet L$, because the additional terms due to the application of such
a Lie series operator are neglected in the secular approximation.

We can finally introduce our secular model up to order $2$ in the
masses by setting
\begin{equation}
  \label{eq:hamsec}
  H^{({\rm sec})}(D_2,\boldsymbol{\xi},\boldsymbol{\eta}) =
  \left\lceil\,
  \langle\,\widetilde{\kern-2pt H}\,\rangle_{\vet{\lambda}}
  \Big|_{\Lbf= \vet{0}}\,\right\rceil_{N_S}\, ,
\end{equation}
i.e., we take the averaged expansion (over the fast angles
$\vet{\lambda}$) of the part of $\,\widetilde{\kern-2pt H}$ that is both
independent from the actions $\Lbf$ and truncated up to a total order
of magnitude $N_S$ in eccentricity and inclination. Since $D_2$ is
$\Oscr\big(e_1^2+i_1^2+e_2^2+i_2^2\big)$, this means that we keep the
Hamiltonian terms $h^{(\Pscr)}_{s;0,j_2}$ with $2s+j_2\le N_S$. From
now on, the parameter $D_2$ is replaced by its explicit value that is
calculated as a function of the initial conditions, so that we can
write the Hamiltonian as follows
\begin{equation}
\label{frm:hsec}
H^{({\rm sec})}(\vet \xi,\vet \eta) = \sum_{s=1}^{N_S/2} h^{({\rm sec})}_{2s}(\vet \xi,\vet \eta),
\end{equation}
where $h_{2s}$ is an homogeneous polynomial of degree $2s$. This means that the expansion contains just
terms of even degree, as a further consequence of the well known D'Alembert rules.

\begin{table}[!]
\begin{center}
\begin{tabular}{l c c}
\hline \hline & $\upsilon$~And~\emph{c} & $\upsilon$~And~\emph{d} \\ 
$m$ [$M_J$] & $15.9792$ & $ 9.9578 $ \\ $a$ [AU] & $0.829$ & $2.53$ \\ $e$ & $0.239$ & $0.31$\\ $i$ [$^\circ$]& $6.865$ & $25.074$\\ 
$M$ [$^\circ$] & $355$ & $335$\\ $\omega$
[$^\circ$] & $245.809$ & $254.302$\\ $\Omega$ [$^\circ$] &
$229.325$ & $7.374$\\ \hline
\end{tabular}
\end{center}
\caption[]{Parameters of the two main planets of
  $\upsilon$~Andromed{\ae} planetary system. The mass of the star
  $\upsilon$~Andromed{\ae} A is equal to $1.31$ solar masses.}
\label{tab:dopo-griglia}
\end{table}

\subsection{Poincar\'e sections for the secular Hamiltonian flow}
\label{sec:Poinc-sect}

We begin to focus on the specific three-body model including the
$\upsilon$~Andromed{\ae} star and the two main exoplanets orbiting
around it. The values for the parameters and the initial conditions
are fixed according to Table~\ref{tab:dopo-griglia}.  These initial
conditions have been selected from the ones compatible with the
observations reported in~\cite{McArt-et-al-2010}, according to a
criterion of robustness. Indeed, the initial conditions are expressed
with large uncertainties, in particular for what concerns the mean
anomalies, that are completely unknown, and the planetary
masses. Therefore, we have chosen between them the ones that minimize
the excursions in eccentricities of the orbits. The complete
methodology will be discussed in a forthcoming paper\footnote{A
  complete explanation of the adopted criterion is included in the
  Ph.D. thesis of C. Caracciolo, see footnote~\ref{phdtesi-Chiara} for
  a more complete reference.}.

We now briefly describe our choices for the truncation parameters that
have been introduced during the discussion of the secular model in
subsection~\ref{sec:secular-model}. There, we have explained the need
of taking into account the main quasi-resonances of the
problem. Concerning the expanded Hamiltonian of the particular case of
$\upsilon$~Andromed{\ae}, the first quasi-resonance that we meet is
$5:1$. Thus, we have chosen to fix the values of the parameters in
such a way that $N_S = 8$ and $K_F = 9$. Considering these settings,
the procedure described in the previous two sub-sections and the
values listed in Table~\ref{tab:dopo-griglia}, the
expansion~\eqref{frm:hsec} is completely defined.

\begin{figure}
\centering
\includegraphics[width=8.0cm]{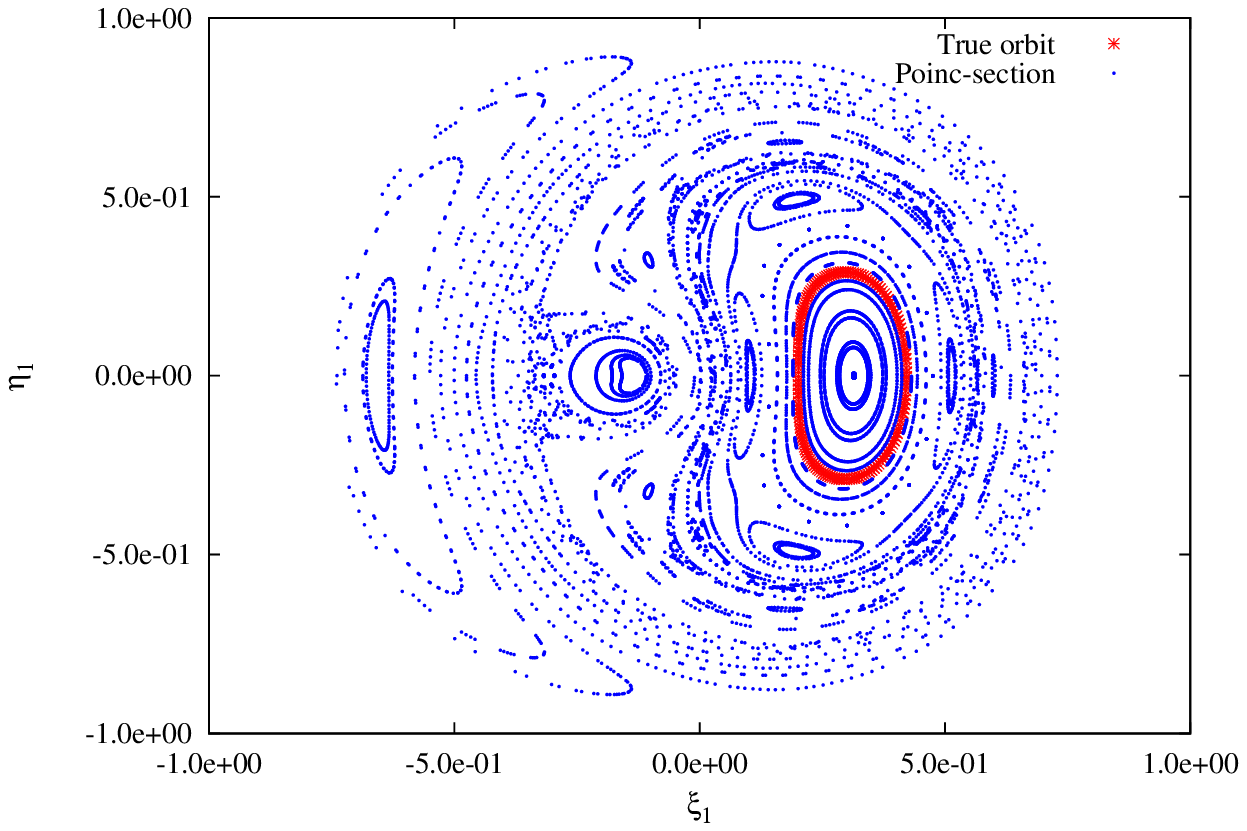}
\includegraphics[width=7.5cm]{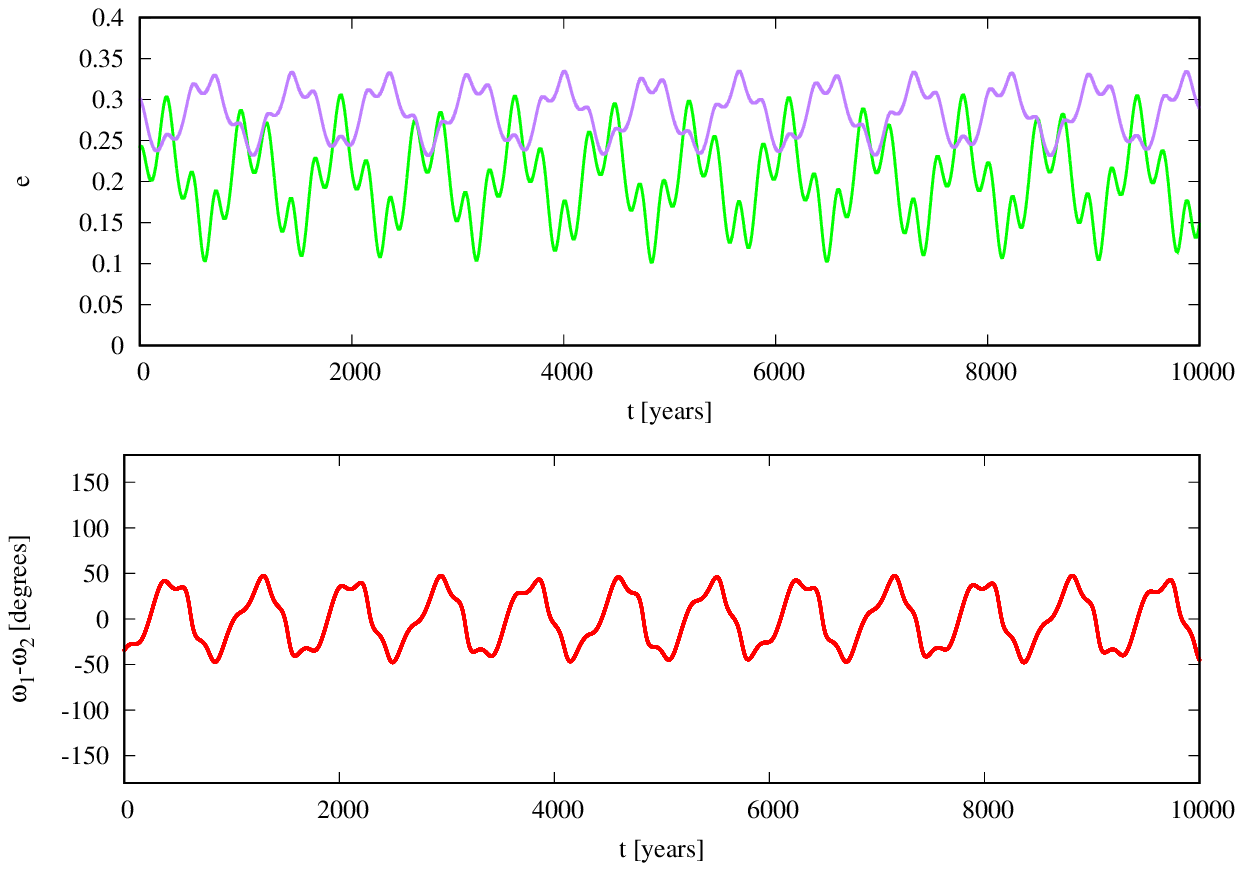}
\caption{Study of the Hamiltonian dynamics defined by $H^{({\rm
      sec})}$ in formula~\eqref{eq:hamsec}. On the right, Poincar\'e
  sections: the orbit with initial conditions compatible with the
  observations is plotted in red, while the orbits referring to other
  initial conditions at the same energy level are in blue. In the
  top-right box, behavior of the planetary eccentricities as a
  function of time: the plots in green and in purple refer to
  $\upsilon$~And~\emph{c} and $\upsilon$~And~\emph{d}, resp. On
  bottom-right, evolution of the difference of the pericenters
  arguments. All the curves reported on the right refer just to the
  {\it true orbit} characterized by the initial values reported in
  Table~\ref{tab:dopo-griglia}.}
\label{fig:poin-sec}
\end{figure}

In the first box of Figure~\ref{fig:poin-sec}, we report the
Poincar\'e sections of the orbital motions induced by the Hamiltonian
flow of $H^{({\rm sec})}$.  They are plotted in correspondence of the
hyperplane $\eta_2=0$, with $\xi_2>0$. They refer to both the initial
values of the orbital elements reported in
Table~\ref{tab:dopo-griglia} (characterizing what we will call,
hereafter, {\it true secular orbit} or, for short, {\it true orbit})
and several other initial conditions having the same energy. Looking
at Figure~\ref{fig:poin-sec}, a few remarks are in order; they are
going to play a fundamental role in the subsequent design of our
approach. Firstly, for most of the Poincar\'e sections the chaotic
effects are not easily remarkable and the two most regular zones are
in the neighborhood of two fixed points (located on the axis of the
abscissas, i.e., with $\eta_1=0$). Moreover, the {\it true secular
  orbit} describes a closed curve that is around one of these points
and does not include the origin. This is in agreement with the fact
that the difference of the arguments of the pericenters librates
without completing a full rotation, as it can be seen in the plot on
the bottom--right of Figure~\ref{fig:poin-sec}, where the dynamical
evolution of the orbital elements is computed by applying, in reverse
order, all the definitions introduced in the previous
subsections~\ref{sbs:3body-model}--\ref{sec:secular-model}.  Finally,
the plots of the eccentricities in the top--right box of
Figure~\ref{fig:poin-sec} highlight that their values oscillate in a
range between $0.1$ and $0.35$.  This makes clear that we have to
develop an approach substantially different with respect to that
described in~\cite{Vol-Loc-San-2018}, which has been designed to study
planetary models similar to that including Sun, Jupiter and Saturn,
where the eccentricities are always smaller than $0.1$ and the
difference of the perihelion arguments is in a regime of full
rotation.

\section{Librational KAM tori for the secular model}
\label{sec:toro}
In this section we describe the technique for the construction of KAM
tori in the case of the secular Hamiltonian model defined in
formula~\eqref{eq:hamsecazang}, with a specific focus on the one which
corresponds to the {\it true secular orbit}, according to the definition
introduced in the previous section.

In order to easily study the model, it is useful to define a set of
action-angle variables $(\vet J, \vet \psi)$ via the canonical
transformation
\begin{equation}
\xi_j=\sqrt{2J_j}\cos \psi_j\ , \qquad \eta_j=\sqrt{2J_j}\sin \psi_j\ ,
\qquad \forall\ j=1,2,
\label{eq:azang}
\end{equation}
being $(\vet \xi, \vet \eta)$ the variables introduced in order to
properly define the secular model. It is important to recall that the
angles $(\psi_1, \psi_2)$ associated to these secular variables are
nearly equal to the arguments of the pericenters $(\omega_1,
\omega_2)$, apart from a small correction due to the (close-the-identity) transformation of coordinates induced by the application of
the Lie series $\exp \Lie{\chi_1^{(\Oscr2)}}$ to the Hamiltonian of
the three-body planetary problem.  After this canonical change of
coordinates, the Hamiltonian~\eqref{frm:hsec} reads
\begin{equation}
  H^{({\rm I})}(\boldsymbol{J},\boldsymbol{\psi})=
\sum_{\ell=1}^{N_S/2} h_{\ell}^{({\rm
      I})}(\boldsymbol{J},\boldsymbol{\psi}) \ ,
\label{eq:hamsecazang}
\end{equation}
where $h_{\ell}$ is an homogeneous polynomial function of degree
$2\ell$ in the square roots of actions $\boldsymbol{J}$ and a
trigonometric polynomial of degree $2\ell$ in angles
$\boldsymbol{\psi}\,$, i.e., it writes
\begin{equation}
  h_{\ell}^{({\rm I})}(\vet{J},\vet{\psi})=
  \sum_{i_1+i_2=2\ell}\,\sum_{j_1=0}^{i_1}\sum_{j_2=0}^{i_2}
  c_{\ell;i_1;i_2;j_1;j_2}^{({\rm I})}
  \sqrt{J_1^{i_1}J_2^{i_2}}\cos\big[(i_1-2j_1)\psi_1+(i_2-2j_2)\psi_2\big]\,.
\end{equation}
The occurrence of only cosines in the formula above with particular
Fourier harmonics in their arguments is again due to the relations
induced by the D'Alembert rules.

It is now convenient to introduce a new set of variables $(\vet I,
\vet \phi)$, where one of the angles is the difference between the
pericenter arguments of the two planets, while the actions are
defined so as to make the change of coordinates canonical; therefore,
the transformation of coordinates is expressed as follows:
\begin{equation}
\phi_1 = \psi_1 -\psi_2\ , \quad \phi_2 = \psi_2\ , \quad
I_1 = J_1\ ,  \quad I_2 = J_2 + J_1\ .
\label{eq:def-azang-Iphi}
\end{equation}
Moreover, let us now introduce the canonical polynomial variables
$(\vet x, \vet y)$ defined as
\begin{equation}
x_j = \sqrt{2 I_j}\cos \phi_j\ , \qquad
y_j=\sqrt{2 I_j}\sin \phi_j\ , \qquad \forall\ j=1,2\ .
\label{eq:def-dalembcoord-xy}
\end{equation}

Let us now remark that making Poincar\'e sections with respect to the
hyperplane $\eta_2=0$, when $\xi_2>0$ is equivalent to impose
$\psi_2=0$, because of the definitions in~\eqref{eq:azang}. Therefore,
looking at
formul{\ae}~\eqref{eq:def-azang-Iphi}--\eqref{eq:def-dalembcoord-xy},
one can easily realize that the drawing on the left of
Figure~\ref{fig:poin-sec} can be seen as a plot of the Poincar\'e
sections in coordinates $(x_1\,,\,y_1)$ with respect to $y_2=0$ and
with the additional condition $x_2>0$. Revisiting the plot of the {\it
  true secular orbit} in the bottom--right box of
Figure~\ref{fig:poin-sec} and in the context of the new canonical
variables is extremely interesting, because it makes clear that
$\phi_1$ is librating around the origin. In fact, we have that $\phi_1
= \psi_1 -\psi_2 \simeq \omega_1 -\omega_2\,$, because the relation
between these differences of angles is given by the
transformation induced by the application of the Lie series $\exp
\Lie{\chi_1^{(\Oscr2)}}$, that is close to the identity. By comparison
with the behavior of the difference between the pericenter arguments
in Figure~\ref{fig:poin-sec}, it is natural to expect that also the
angle $\phi_1$ is librating in a range that cannot be much larger than
$[-50^\circ,50^\circ]$.

Focusing on the axis of the abscissas of the Poincar\'e sections
plotted in blue in Figure~\ref{fig:poin-sec}, one can remark that
there is a fixed point located at $(\xi_1\simeq 0.3\,,\,\eta_1=0)$;
moreover, it is surrounded by orbits lying on nearly elliptical closed
curves. It can be seen as a periodic orbit, namely a lower dimensional
(here, one-dimensional) torus that in this case is transversely
elliptic. From an astronomical point of view, it is characterized by
the fact that the pericenters stay nearly perfectly
anti-aligned\footnote{These angles are expressed in the Laplace
  reference frame, where the longitudes of the nodes are always
  opposite. This is the reason why the alignment of the arguments of
  the pericenters actually corresponds to the anti-alignment of the
  pericenters and {\it vice versa}.}, being
$0=\phi_1\simeq\omega_1-\omega_2\,$.  The {\it true secular orbit} is
in libration around such an extremely peculiar configuration. The
Poincar\'e sections in Figure~\ref{fig:poin-sec} are structured in
such a way that it is natural to expect that the periodic solution
corresponding to the anti-alignment of the pericenters is extremely
robust. This property should allow the one-dimensional elliptic torus
to influence the dynamics and the stability in its neighborhood and,
in particular, for what concerns the {\it true secular
  orbit}. Therefore, we will proceed by first constructing that
specific periodic solution which describes a full rotation of the
angle $\phi_2\,$, without any oscillation in the transverse directions
corresponding to the pair of canonical variables $(I_1\,,\,\phi_1)$
or, equivalently, $(x_1\,,\,y_1)$.

\subsection{Construction of the lower dimensional elliptic torus}
\label{subsec:ell}
The construction of the wanted periodic orbit can be done with a
normal form algorithm\footnote{In~\cite{Caracciolo-2021}, the
  convergence of such a constructive algorithm is proved under
  suitable hypotheses.} similar to the one described
in~\cite{San-Loc-Gio-2011} and~\cite{Car-Loc-2021}, after having
defined suitable transformations of coordinates in order to be
approximately centered around the solution we aim to build. The
location of the wanted periodic orbit can be determined numerically
with a method that is based on the frequency analysis~(see,
e.g.,~\cite{Laskar-03}), in short FA, hereafter. In the following, we
briefly recall the main features of such a computational procedure,
that is accurately described, e.g., in~\cite{Car-Loc-2021}.

We look for a set of initial conditions close to an orbit in which the
angle $\phi_2$ rotates periodically, without oscillations in the
transverse plan $(x_1, y_1)$. The {\it true orbit} is near that
solution as there is a small libration of the angle
$\phi_1$. Therefore, we start by extracting from the numerical signal
$x_2(t) + \imunit y_2(t)$, which is related to the motion of the angle
$\phi_2$ of the {\it true orbit}, its main frequency $\nu_{2,1}\,$,
being $\sum_{s=1}^{\Nscr_C} A_{2,s} e^{\imunit (\nu_{2,s}t +
  \phi_{2,s})}$ the decomposition obtained by using FA, where
$\Nscr_C$ is the number of components considered and, hereafter, we
assume that the amplitudes $A_{j,s} > 0$ are in descending order with
respect to $s=1, \ldots, \Nscr_C \ \forall \ j=1,2$. Then, we focus on
the signal of the other pair of canonical coordinates, i.e., $x_1(t) +
\imunit y_1(t) \simeq \sum_{s=1}^{\Nscr_C} A_{1,s} e^{\imunit
  (\nu_{1,s} t+ \phi_{1,s})}$, that is related to the librational
motion of the angle $\phi_1$; therefore, we expect to find
contributions given by frequencies that are not integer multiples of
$\nu_{2,1}\,$. We define as a new initial condition $\tilde x_j(0) +
\imunit \tilde y_j(0) = \sum_{s=1}^{\Nscr_C}\tilde A_{j,s} e^{\imunit
  \phi_{j,s}}$, where we are considering only the components $\tilde
A_{j,s}$ related to the frequencies $\nu_{j,s}$ that are integer
multiples of $\nu_{2,1}$, otherwise we put the corresponding amplitude
$A_{j,s} = 0$. We integrate the equation of motion with this other
initial condition and we iterate the procedure until both the signals
can be decomposed, up to a tolerance threshold, in terms of a single
frequency $\nu_{2,1}$.

Once that we have found the initial condition $( \vet x^\star, \vet
y^\star)$ for a periodic solution, we are able to introduce a first
translation on the action $I_2$ by defining $p= I_2 -I^\star$, where
$I^\star = (x_2^{\star 2}+y_2^{\star 2})/2$, and to expand the
Hamiltonian in Taylor series around $p=0$. Then, we divide the
variables in two different couples: we rename the angle $\phi_2$ as
$q$, so that $(p,q)\in \reali\times \toro$ is the action-angle couple
describing the periodic motion, while we still use the polynomial
variables $(x_1,y_1)$ for the motion transverse to the periodic
orbit. The last preliminary translation is on $x_1\,$, in order to
have expansions around the value $x_1^\star$, given by the initial
condition computed numerically.  Let us emphasize that, since the
fixed point we are trying to approximate in Figure~\ref{fig:poin-sec}
corresponds to $\phi_1=0$, we have that $y_1= 0$ and here a translation
is not needed. Finally, before starting with the construction of a
better approximation of the elliptic torus by putting the Hamiltonian
in a suitable normal form, we rescale the transverse variables $(\bar
x_1,y_1)$, being $\bar x_1= x_1 - x_1^\star$, in such a way that the
quadratic part in the new variables $(x,y)$ is in the form
$\Omega^\star(x^2 + y^2)/2$. This rescaling can be done by a canonical
transformation as the quadratic part does not have any mixed term
$\bar x_1 y_1$ and the coefficients of $\bar x_1^2$ and $y_1^2$ have
the same sign, because of the proximity to an elliptic equilibrium
point. Thus, since such a quadratic part is in the preliminary form
$a\bar x_1^2+ b y_1^2$, it suffices to define the new variables
$(x,y)$ as $x = \sqrt[4]{\frac a b}\, \bar x_1, \ y = \sqrt[4]{\frac b
  a}\, y_1$.

After all these transformations, the expansion of the Hamiltonian can
be written as follows:
\begin{equation}
  \label{frm:ham-ell}
  \begin{aligned}
    H(p,q,x,y) = & E+ \omega^\star \ p +
    \Omega^\star\frac{x^2 +y^2}{2}+
\Fscr_{\rm{h.o.t.}}(p,q,x,y)    
   + \\ &+  \sum_{s> 0}
    f_0^{(s)}(q)+ \sum_{s> 0} f_1^{(s)}(q,x,y) +
    \sum_{s> 0} f_2^{(s)}(p,q,x,y)\, ,
  \end{aligned}
\end{equation}
where each $f_\ell^{(s)}$ is a function in the variables $(p,q,x,y)$,
whose corresponding Taylor-Fourier series is given by a finite number
of terms. More precisely, $\ell$ is its total degree in the square
root of the actions, i.e., $\ell = 2j_1+j_2$, where $j_1$ is the
degree in $p$ and $j_2$ the degree in $(x,y)$, while $2s$ is the
maximum trigonometrical degree in the angle
$q$. In~\eqref{frm:ham-ell}, we have highlighted the terms with total
degree in the square root of the actions up to $2$ and we have
collected all the terms of higher order in $\Fscr_{\rm{h.o.t.}}$. This
choice is coherent with our purpose to start the iteration of an
algorithm for the construction of a lower dimensional elliptic torus,
by following the approach described in~\cite{San-Loc-Gio-2011} and
in~\cite{Car-Loc-2021}; this is feasible provided that the terms
appearing in the second row of~\eqref{frm:ham-ell} are small enough
with respect to those in the first row. Furthermore, $E$ is a constant
term, meaning just the energy level of the invariant manifold $p = x=
y = 0$ in the limit case with $f_0^{(s)} = f_1^{(s)} = f_2^{(s)} = 0
\ \forall\ s > 0$. In the following, we will just recall the
procedure, by emphasizing what is necessary to adapt in such a way to
fit with the present context.

The expansion of the Hamiltonian in~\eqref{frm:ham-ell} can be
visually rearranged as
\begin{equation}
\label{eq:hamell}
  \vcenter{ \halign to\hsize{ &\hfil$\displaystyle\>{#}$\hfil
      &\ \ \hfil$\displaystyle\>{#}$\hfil
      &\ \ \hfil$\displaystyle\>{#}$\hfil
      &\ \ \hfil$\displaystyle\>{#}$\hfil
      &\ \ \hfil$\displaystyle\>{#}$\hfil
      &\ \ \hfil$\displaystyle\>{#}$\hfil
      &\ \ \hfil$\displaystyle\>{#}$\hfil \cr & &\vdots
      &\vdots&\vdots&\vdots&\vdots \cr & &f_3^{(0,0)}(p,x,y)
      &f_3^{(0,1)}(p,q,x,y)  &\ldots &f_3^{(0,s)}(p,q,x,y)&\ldots \cr & & \omega^{(0)} p + \Omega^{(0)} \frac{x^2
        + y^2}{2} &f_2^{(0,1)}(p,q,x,y)  &\ldots
      &f_2^{(0,s)}(p,q,x,y) &\ldots \cr &H^{(0)}\,=\, \Bigg.\sum & 0 &f_1^{(0,1)}(q,x,y)  &\ldots
      &f_1^{(0,s)}(q,x,y)  &\ldots &\ , \cr & &E^{(0)}\,
      &f_0^{(0,1)}(q)\, &\ldots &f_0^{(0,s)}(q)\, &\ldots \cr } }
\end{equation}
where $\omega^{(0)}= \omega^*$, $\Omega^{(0)}=\Omega^*$ and
$E^{(0)}=E$, the first upper index of the polynomial terms and of the
Hamiltonian refers to the normalization step, while the second one
refers to the trigonometrical degree in $q$. In the normalization
procedure, we will define a sequence of canonical transformations,
with the aim of removing the contributions appearing in the first
three rows, that are terms having total degree in the square root of
the actions up to $2$, except for the part in normal form that does
not depend on the angle $q$. In order to avoid the proliferation of
symbols, by abuse of notation we are going to use the same name for
the variables, even if we are indeed applying several changes of
coordinates.

The first change of coordinates is identified by the generating
function $\chi_1^{(1)}(q)$, of trigonometrical degree up to $2$ in the
angle $q$, which solves the following homological equation:
\begin{equation}
\label{frm:1eqhom}
\left\{\chi_0^{(1)},\ \omega^{(0)} p\right\} + f_0^{(0,1)}(q) =
\langle f_0^{(0,1)}(q)\rangle_q.
\end{equation}
Let us stress that the term in the r.h.s. of the equation above is
constant because it denotes the angular average of $f_0^{(0,1)}$;
therefore, it can change the energy value but it has not any
role in the Hamilton equations.

We will obtain an Hamiltonian $ H^{(I; 1)}= \exp \lie{\chi_0^{(1)}}
H^{(0)}$ with $f_0^{(I;1,1)}(q)=0$, where the new upper index $I$ is
now referred to the first substep and $E^{(1)} = E^{(0)}+ \langle
f_0^{(0,1)}(q)\rangle_q$.  In order to better understand the
following, it is convenient to imagine the expansion of $H^{(I;1)}$ as
in formula~\eqref{eq:hamell}, by replacing each term $f_\ell^{(0,s)}$
with $f_\ell^{(I;1,s)}$, having the same functional properties, i.e.,
$\ell$ is the total degree in the square root of the actions, while
$2s$ is the maximum trigonometrical degree in the angle $q$.

The second substep is meant to eliminate terms linear in the square
root of the actions. We introduce $\chi_1^{(1)}(q,x,y) $, linear in
$(x,y)$ and of trigonometrical degree up to $2$ in the angle $q$,
such that
\begin{equation}
\label{frm:2eqhom}
\left\{\chi_1^{(1)}, \ \omega^{(0)} p+ \Omega^{(0)}
\frac{x^2+y^2}{2}\right\} + f_1^{(I;1,1)}(q,x,y)  = 0.
\end{equation}
In this case, since the origin of the transverse variables corresponds
to an elliptic equilibrium point, terms that are linear in $x$ and $y$
with non-zero angular average cannot appear in the r.h.s. of the
equation above. Therefore, we are able to define $H^{(II;1)} = \exp
\lie{\chi_1^{(1)}} H^{(I;1)}$ with $f_1^{(II;1,1)}(q,x,y) =0$. In
order to better understand the following, once again, it is convenient
to imagine the expansion of $H^{(II;1)}$ as in
formula~\eqref{eq:hamell}, by replacing each term $f_\ell^{(0,s)}$
with the corresponding one, that is denoted with $f_\ell^{(II;1,s)}$
and has the same functional properties.

The third substep aims to remove terms linear in the actions that
depend on the angle $q$ and it is split in two different stages, each
of them responsible for the removal of terms linear in $p$ or
quadratic in $(x,y)$. We define $X_2^{(1)}(p,q)$, linear in $p$ and of
trigonometrical degree $2$ in the angle $q$, such that
\begin{equation}
\label{frm:3eqhom1}
\left\{X_2^{(1)}, \ \omega^{(0)} p\right\} + f_2^{(II;1,1)}(p,q) = \langle
f_2^{(II;1,1)}(p,q)\rangle_q;
\end{equation}
we define $Y_2^{(1)}(q,x,y) $, quadratic in $(x,y)$ and of
trigonometrical degree up to $2$ in the angles $q$, such that
\begin{equation}
\label{frm:3eqhom2}
\left\{Y_2^{(1)}, \ \omega^{(0)} p+ \Omega^{(0)}
\frac{x^2+y^2}{2}\right\} + f_2^{(II;1,1)}(q,x,y)   = \langle
f_2^{(II;1,1)}(q,x,y) \rangle_q.
\end{equation}
In both the previous formul{\ae} we have denoted the kind of terms by
writing explicitly the dependence on the variables.  The terms that
we cannot remove with these transformations have non-zero average on
$q$; moreover, they are linear in the actions and essentially of the
same type of the normal form. Therefore, they are added to the normal
form part itself by defining new frequencies $\omega^{(1)}$ and
$\Omega^{(1)}$ as follows:
\begin{align}
\omega^{(1)}p &=\omega^{(0)}p + \langle f_2^{({\rm II};1,1)}(p,q)\rangle_{q}\,,\\
\Omega^{(1)}\frac{x^2 +y^2}{2} &=\Omega^{(0)}\frac{x^2 +y^2}{2}+ \langle f_2^{({\rm II};1,1)}(q,x,y) \rangle_{q}\,,
\label{frm:trasv-freq}
\end{align}
where in the last formula, $\langle f_2^{({\rm II};1,1)}(q,x,y)
\rangle_{q}$ is put in diagonal form. This can be done perturbatively
by using the Lie transform, i.e., by performing an infinite sequence
of Lie series (see~\cite{Gio-Loc-San-2014} and~\cite{Car-Loc-2021} for
further details). Therefore, the first normalization step is completed
by applying $\exp \lie{X_2^{(1)}}$, $\exp \lie{Y_2^{(1)}}$ and the Lie
transform (that has been mentioned just above) to the intermediate
Hamiltonian $H^{(II;1)}$. This allows us to fully determine $H^{(1)}$,
whose expansion can be represented as in formula~\eqref{eq:hamell}, by
replacing each term $f_\ell^{(0,s)}$ with the corresponding one, i.e.,
$f_\ell^{(1,s)}$, having the same functional properties. Of course,
for what concerns the new Hamiltonian $H^{(1)}$, in the expansion
analogous to that written in~\eqref{eq:hamell}, the angular velocities
$\omega^{(1)}$, $\Omega^{(1)}$ and the the energy level $E^{(1)}$ will
appear in place of $\omega^{(0)}$, $\Omega^{(0)}$ and $E^{(0)}$,
respectively; moreover, in the first three cells of the second row,
$f_\ell^{(0,1)}$ shall be replaced by $f_\ell^{(1,1)}=0$ for
$\ell=0,\,1,\,2$, accordingly with the definitions of the homological
equations~\eqref{frm:1eqhom}--\eqref{frm:3eqhom2}.

The generic $r$-th step is performed in the same way, with the only
difference that the generating functions $\chi_0^{(r)}$,
$\chi_1^{(r)}$, $X_2^{(r)}$ and $Y_2^{(r)}$ are determined in such a
way to remove the perturbing terms of trigonometrical degree up to
$2r$. Let us stress the fact that, in order to solve the homological
equations and to let the whole procedure to be convergent, some
non-resonance conditions have to be satisfied. These requirements can be
summarized by the following formula:
\begin{equation}
\min_{ {0 \le | k | \le 2r , \ 0 \le | \ell| \le 2} \atop { | k | +
    |\ell|> 0 } } \Big | k \omega^{(r-1)}+ \ell \Omega^{(r-1)} \Big |
\ge \frac{\gamma}{| k|^{\tau} +1} \ \ { \rm with} \ \gamma > 0, \ \tau
> 0 \ ,
\end{equation}
where let us emphasize that here, for $\ell =0$, we do not have small
divisors, since we deal with a single frequency (that is related to
the wanted periodic orbit) instead of a vector.

The convenience of such a procedure is highlighted by the fact that,
after an infinite number of iterations of the algorithm, the
Hamiltonian reads
\begin{equation}
  \label{frm:ham-inf-elltor}
  H^{(\infty)} = E^{(\infty)} +\omega^{(\infty)}\ p +
  \Omega^{(\infty)}\, \frac{x^2 + y^2}{2}+ \sum_{s\ge
    0} \sum_{\ell > 2}f_\ell^{(\infty, s)}(p,q,x,y)
\end{equation}
and, therefore, it is evident that $(0, q_0+\omega^{(\infty)}t, 0 ,0)$
is a periodic orbit. Besides, from a practical point of view, we
iterate the algorithm only for a finite number of steps and we
construct an approximation of the periodic orbit, being the final
Hamiltonian after $\bar r$ steps
\begin{equation}
  \label{frm:ham-barr-elltor}
  H^{(\bar r)} = E^{(\bar r)} +\omega^{(\bar r)}\ p +
  \Omega^{(\bar r)}\, \frac{x^2 + y^2}{2}+ \sum_{s\ge
    0} \sum_{\ell > 2}f_\ell^{(\bar r,s)}(p,q,x,y) + \sum_{s>
    \bar r} \sum_{\ell \ge 0}^2f_\ell^{(\bar r,s)}(p,q,x,y)\ ,
\end{equation}
where $\sum_{s> \bar r} \sum_{\ell \ge 0}^2f_\ell^{(\bar r,
  s)}(p,q,x,y) $ is the remainder of the normal form.

Let us conclude this subsection with a comment about the choice of the
initial translations $I^\star$ and $x_1^\star$. Indeed, if they are
accurate enough, we expect that the algorithm will converge to an approximated
periodic orbit defined by $(\omega^{(\bar r)},\ \Omega^{(\bar r)})$
close to the one we are aiming at, namely the one we have found by
using the FA. However, at this level, the value of $I^\star$ can be
improved, in order to ``automatically" lead to a periodic orbit as close as possible to the numerical one. This is done by using a
Newton method to solve the equation $E(I) = \bar E$, where $E(I)$ is
the energy associated to the (approximated) periodic orbit we have
found by using the translation $I$ and iterating the algorithm for a
finite number of steps, while $\bar E$ is the energy of the periodic
orbit found by using the FA. The same Newton method can be used to
improve the translation $x_1^\star$ and, therefore, to remove from the
Hamiltonian the terms linear in $(\bar x_1,y_1)$: the equation to be
solved in this case is $\frac{\partial \langle
  H\rangle_q}{\partial{x}} (0,0,0,0) = 0$, where we recall that $x =
\sqrt[4]{\frac a b}\, \bar x_1$. The initial translations $I^\star$
and $x_1^\star$ as defined at the beginning of this subsection are in
general good enough initial guesses for ensuring the simultaneous
convergence of both the Newton methods.

\subsection{Construction of the librational KAM torus}
\label{sec:kamtor}
The normal form related to the elliptic torus is now considered as the
starting point for the construction of the full dimensional KAM tori
around it and, in particular, the invariant surface corresponding to
initial conditions compatible with the observations. Since we cannot
derive the secular frequencies by the observations, we use again the
FA to compute the frequency vector $\tilde{\vet \omega}$, being
$\tilde \omega_1$ and $\tilde \omega_2$ the frequencies associated to
the rotation of the angle $\phi_2$ and to the oscillations of the
angle $\phi_1$.

Initially, we translate again the coordinates, in order to be centered
around a good approximation of the \emph{true orbit}. The variables
$(x,y)$ for which the Hamiltonian is in normal form for elliptic tori
are centered around the equilibrium point, but we can identify the
right translation on the action $I = (x^{2}+y^{2})/2$ by expressing
the initial conditions of the {\it true orbit} in these
coordinates\footnote{The algorithm described in the previous
  subsection provides the explicit change of coordinates needed to
  introduce the normal form and its inverse. A more detailed
  description of how to define them is deferred to the next section.};
then, we can expand the Hamiltonian in Taylor series around $p_2=0$,
being $p_2 = I-I_2^\star$, where $I_2^\star$ is the value of the
action at time $0$ of the {\it true orbit}. Moreover, we rename the
angle associated to the transverse polynomial variables $(x,y)$ as
$q_2$ and we recollect the variables with the other action-angle
couple $(p,q)$ (renamed as $(p_1,q_1)$), so that the Hamiltonian is
now expressed in the variables $(\vet p, \vet q)$, being the first
couple referred to the periodic motion of the angle $\phi_2$, while
the second one is referred to the libration as it is measured from the
elliptic torus.

We are ready to proceed with a classic algorithm for the construction
of the Kolmogorov normal form, by following the approach described in
section~4 of~\cite{Gab-Jor-Loc-2005}. First, the expansion of the
initial Hamiltonian $\Hscr^{(0)}$ (that has been obtained from
$H^{(\bar r)}$ written in~\eqref{frm:ham-barr-elltor}, by applying the
translation described just above) can be visually reorganized as
follows:
\begin{equation}
\label{eq:hamexp}
  \vcenter{ \halign to\hsize{ &\hfil$\displaystyle\>{#}$\hfil
      &\ \ \hfil$\displaystyle\>{#}$\hfil
      &\ \ \hfil$\displaystyle\>{#}$\hfil
      &\ \ \hfil$\displaystyle\>{#}$\hfil
      &\ \ \hfil$\displaystyle\>{#}$\hfil
      &\ \ \hfil$\displaystyle\>{#}$\hfil
      &\ \ \hfil$\displaystyle\>{#}$\hfil \cr & &\vdots
      &\vdots&\vdots&\vdots&\vdots \cr & &f_2^{(0,0)}(\boldsymbol{p})
      &f_2^{(0,1)}(\boldsymbol{p},\boldsymbol{q}) &\ldots
      &f_2^{(0,s)}(\boldsymbol{p},\boldsymbol{q}) &\ldots \cr
      &\Hscr^{(0)}(\boldsymbol{p},\boldsymbol{q})\,=\, \Bigg.\sum
      &\boldsymbol{\omega}^{(0)}\cdot\boldsymbol{p}
      &f_1^{(0,1)}(\boldsymbol{p},\boldsymbol{q}) &\ldots
      &f_1^{(0,s)}(\boldsymbol{p},\boldsymbol{q}) &\ldots &\qquad ,
      \cr & &E^{(0)}\, &f_0^{(0,1)}(\boldsymbol{q})\, &\ldots
      &f_0^{(0,s)}(\boldsymbol{q})\, &\ldots  \cr}}
\end{equation}
being the generic term $f_j^{(0,s)}$ an homogeneous polynomial of
degree $j$ in the actions $\boldsymbol{p}$ and of trigonometrical degree
up to $2s$ in $\boldsymbol{q}$.  The Kolmogorov's normalization
algorithm requires to remove all the terms of the
Hamiltonian~\eqref{eq:hamexp} of degree~$0$ or~$1$ in the actions
$\boldsymbol{p}$, with the exception of the term
${\boldsymbol{\omega}^{(0)}}\cdot\boldsymbol{p}\,$. Let us stress that
an accurate definition of the translation vector $I_2^\star$ should
lead to a frequency vector $\vet{\omega}^{(0)} \simeq \tilde{\vet
  \omega}$.

In order to construct the normal form, we start by determining the
generating function $X^{(1)}(\vet q)$ such that
\begin{equation}
\label{eq:kolgen1}
  \poisson{X^{(1)}}{\boldsymbol{\omega}^{(0)}\cdot\boldsymbol{p}}+
  f_0^{(0,1)}=\langle f_0^{(0,1)} \rangle_{\vet q}\ ,
\end{equation}
where $X^{(1)}$ is a trigonometric polynomial of degree~$2$. The term
$\langle f_0^{(0,1)} \rangle_{\vet q}$ will contribute to the energy
as a new constant term, to be summed up to the initial value of the
energy $E^{(0)}$ (again, close to the energy of the {\it true orbit}). We
will then obtain a new Hamiltonian $\hat
\Hscr^{(1)}=\exp\Lscr_{X^{(1)}}\Hscr^{(0)}$, whose expansion can be
visualized as in~\eqref{eq:hamexp} with generic terms now denoted by
$\hat f_j^{(1,s)}$, being $\hat f_0^{(1,1)}=0$ as a consequence of
equation~\eqref{eq:kolgen1}.

We proceed in an analogous way to complete this first Kolmogorov's
normalization step: we compute the generating function
$\chi_2^{(1)}(\boldsymbol{p},\boldsymbol{q})$ such that
\begin{equation}
\label{eq:kolgen2}
  \poisson{\chi_2^{(1)}}{\boldsymbol{\omega}^{(0)}\cdot\boldsymbol{p}}+
  \hat f_1^{(1,1)}=\langle\hat f_1^{(1,1)}\rangle_{\boldsymbol{q}}\ ;
\end{equation}
therefore, $\chi_2^{(1)}$ will be linear in $\boldsymbol{p}$ and of
trigonometrical degree equal to~2 in $\boldsymbol{q}$. Let us remark
that it is possible to solve the previous homological
equations~\eqref{eq:kolgen1} and~\eqref{eq:kolgen2}, provided that
$|\vet{k} \cdot \vet{\omega}^{(0)}| > 0$ for $\vet{k}\in\interi^2$
with $|\vet{k}|=1,2$\,, being $|\vet{k}| = |k_1| + |k_2|$\,.

If we do not introduce a further change of coordinates with the aim of
fixing the frequency vector of the torus in construction, we will need
to define a new frequency vector $\boldsymbol{\vet\omega}^{(1)}$ in
such a way that
\begin{equation}
\boldsymbol{\omega}^{(1)}\cdot\boldsymbol{p}=%
\boldsymbol{\omega}^{(0)}\cdot\boldsymbol{p} +%
 \langle\hat f_1^{(1,1)}\rangle_{\boldsymbol{q}}\ .
\end{equation}
Therefore, we will obtain the new Hamiltonian
$\Hscr^{(1)}=\exp\Lscr_{\chi_2^{(1)}}\hat \Hscr^{(1)}$, whose
expansion can be again visualized as in~\eqref{eq:hamexp}, where each
generic term $f_j^{(0,s)}$ is replaced with the corresponding one,
i.e., $f_j^{(1,s)}$, having the same functional properties; moreover,
$f_0^{(1,1)}$ and $f_1^{(1,1)}$ are now equal to 0.

The generic $r$-th normalization step can be performed in the same
way, because the Hamiltonian expansion $\Hscr^{(r-1)}$, if we replace
the upper index $0$ with $r-1$, is of the same form as that described
in~\eqref{eq:hamexp}. Although the procedure is rather standard, we
retain useful to report here some additional detail about the
construction, being also necessary for the following discussions about
the convergence of the procedure. The generating functions $X^{(r)}$
and $\chi_2^{(r)}$ are introduced by solving the homological equations
obtained by replacing the upper indices $0$ and $1$ with $r-1$ and
$r$, respectively, in formul{\ae}~\eqref{eq:kolgen1}
and~\eqref{eq:kolgen2} and they can be written as follows:
\begin{equation}
\label{frm:X}
X^{(r)}(\vet q) = \sum_{0<|\vet k|\le 2r} - c_{\vet k}^{(r-1,r)} \frac
{e^{i \vet k \cdot \vet q}}{i\vet k \cdot \vet \omega^{(r-1)}}\ ,
\end{equation}
being $f_0^{(r-1,r)}(\vet q) = \sum_{0<|\vet k|\le 2r} c_{\vet
  k}^{(r-1,r)} e^{i \vet k \cdot \vet q}$ the term that does not
depend on the actions $\vet p $ to be removed at the $r$-th step,
while the second generating function is defined as
\begin{equation}
\label{frm:chi2}
\chi_2^{(r)}(\vet p, \vet q) = \sum_{j=1}^2\sum_{0<|\vet k|\le 2r} -
c_{j, \vet k}^{(r-1,r)} p_j \ \frac {e^{i \vet k \cdot \vet q}}{i\vet k
  \cdot \vet \omega^{(r-1)}}\ ,
\end{equation}
being $\hat f_1^{(r-1,r)}(\vet p, \vet q) = \sum_{j=1}^2 \sum_{0<|\vet
  k|\le 2r} c_{j, \vet k}^{(r-1,r)} p_j \ e^{i \vet k \cdot \vet q} $
the expansion of the unwanted term that is met at the $r$-th step and
is linear in the actions $\vet p$.  Therefore, the homological
equations can be solved provided that the following non-resonance
condition holds true:
\begin{equation}
\label{eq:res}
|\vet{k} \cdot \vet{\omega}^{(r-1)}| > 0\ , \qquad
\forall\ \vet{k}\in\interi^2\setminus\{\vet{0}\}\ {\rm with}
\ |\vet{k}|\le 2r\ .
\end{equation}

\noindent
The new Hamiltonian is finally given by
\begin{equation}
\Hscr^{(r)}= \exp\Lscr_{\chi_2^{(r)}}\hat
\Hscr^{(r)}\quad {\rm with}\quad \hat \Hscr^{(r)}=\exp\Lscr_{X^{(r)}}
\Hscr^{(r-1)}\ ,
\end{equation}
where the effect on the expansions that is due to the application of
the Lie series can be fully described in the following way. In the
spirit of a programming language, we define the intermediate functions
$\hat f_{\ell}^{(r,s)} = f_{\ell}^{(r-1, s)}\ \forall \ \ell, s \ge 0$
and then the contributions due to the Lie derivatives with respect to
the generating function $X^{(r)}$, that is an operator decreasing by
$1$ the degree in the actions and increasing by $2r$ the
trigonometrical degree in the angles, are added in such a way that
\begin{equation}
\label{frm:dopochi1}
\hat f_{\ell -j}^{(r, jr+s)} \hookleftarrow
\frac{1}{j!}\Lie{X^{(r)}}^j f_\ell^{(r-1,s)}\,,
\quad {\rm for}  \ 1 \le j \le \ell,\ s\ge 0\ ,
\end{equation}
where with the notation $a \hookleftarrow b$ we mean that the quantity
$a$ is redefined so as to be equal $a+b$. For what concerns the Lie
derivative with respect to the generating function $\chi_2^{(r)}$, it
does not change the degree in the actions while it increases by $2r$
the trigonometrical degree in the angles, therefore, we can define
$f_{\ell}^{(r, s)}=\hat f_{\ell}^{(r, s)}$ $\forall \ \ell, s \ge 0$
and then add all the other generated terms as follows:
\begin{equation}
\label{frm:dopochi2}
f_{\ell}^{(r, jr+s)} \hookleftarrow
\frac{1}{j!}\Lie{\chi_2^{(r)}}^j \hat f_\ell^{(r-1,s)}\,,
\quad {\rm for}\ j \ge 1, \ \ell \ge 0, \ s\ge 0\ .
\end{equation}

The algorithm is iterated up to a finite number of steps $\bar r$, thus
at the end we have an approximated torus with frequency vector $\vet
\omega^{(\bar r)}$, being the $r$-th Hamiltonian in the form
\begin{equation}
\Hscr^{(\bar r)} = E^{(\bar r)} + \vet \omega^{(\bar r)}\cdot \vet p + \sum_{s \ge 0} \sum_{\ell \ge 2} f_{\ell}^{(\bar r,s)} (\vet p, \vet q)+ \sum_{s> \bar r} \sum_{\ell = 0}^{1} f_{\ell}^{(\bar r,s)} (\vet p, \vet q)\ .
\label{frm:hbarr}
\end{equation}
Once again, we recall that if the initial translation $I_2^\star$ is
accurate enough, we should end up with a torus close to the {\it real}
one, with $\vet \omega^{(\bar r)} \simeq \tilde{\vet \omega}$. Also in
this case, the shift on the action value $I_2^\star$ can be improved
by means of a Newton method; furthermore, it could be useful to
add an appropriate extra translation to the action $p_1$, that, we
recall, is already expanded around $I^\star$, i.e., the initial value
of the action of the periodic orbit. The aim is to find the
translations $\tilde I_1$ and $\tilde I_2$ such that $\vet
\omega^{(\bar r)}(\tilde I_1,\tilde I_2) = \tilde{\vet \omega}$, being
the translation applied before the iteration of the Kolmogorov's
algorithm defined as
\begin{equation}
\label{solo1trasl}
\chi(\vet q) = \tilde{\vet  I}\cdot \vet q \,.
\end{equation}
The initial translations are iteratively adjusted by applying the
following refinement formula:
\begin{equation}
  \tilde{ \vet I}^{(n)} =
  \tilde{\vet I}^{(n-1)} -
  \frac{\Delta \vet \omega(\tilde{\vet I}^{(n-1)})}{\Jscr(\tilde{\vet I}^{(n-1)})}
\end{equation}
where $\Delta \vet \omega(\tilde{\vet I}) =\vet \omega^{(\bar
  r)}(\tilde{\vet I}) - \tilde{\vet \omega}$, $ \tilde{\vet I}^{(0)} =
(0, I_2^\star)$ and the Jacobian $\Jscr(\tilde{\vet I})$ of the
function $\vet \omega(\tilde {\vet I})$ is numerically computed, being
the evaluation of such a function rather implicit.

After having found $\tilde I_1$ and $\tilde I_2$ such that the
approximation of the torus is good enough, we iterate one last time
the algorithm constructing the Kolmogorov normal form
$\Kscr^{(\infty)}$ (written in~\eqref{eq:haminfty}), in this case by
keeping fixed the value of the angular velocity vector so as to be
equal to $\tilde{\vet \omega}$ by means of a translation to be
performed at each normalization step. The algorithm is restarted from
the new initial Hamiltonian $\Kscr^{(0)}$ that is defined as
\begin{equation}
\Kscr^{(0)} = E^{(0)} + \tilde{\vet \omega}\cdot \vet p + \sum_{s \ge 0} \sum_{\ell \ge 2} f_{\ell}^{(0,s)} (\vet p, \vet q)+ \sum_{s\ge 1} \sum_{\ell = 0}^{1} f_{\ell}^{(0,s)} (\vet p, \vet q),
\label{frm:K0}
\end{equation}
where the terms appearing in such a formula are obtained by the
corresponding ones in the expansion of $\Hscr^{(\bar r)}$
in~\eqref{frm:hbarr}, after having translated the actions by a vector
$\vet \xi^{(0)}$ such that $\vet \omega^{(\bar r)} + C^{(0)} \vet
\xi^{(0)} = \tilde{\vet\omega}$, where $ C^{(0)} = {\rm Hess}
\sum_{s=0}^{\bar r} \langle f_2^{(\bar r,s)}\rangle_{\vet q}$. Once
again, we do not change the symbols $E$ and $f$ after having performed
a translation (that is a canonical transformation) by abuse of
notation.  Therefore, the main difference with the algorithm
previously described is that we have to compute another generating
function $\vet \xi^{(r)} \cdot \vet q$ at each step $r$, so as to
solve the following homological equation
\begin{equation}
\label{frm:trasl}
  \left\{\vet \xi^{(r)}\cdot \vet q , \frac 1 2 C^{(r)} \vet p \cdot \vet p \right\} + \vet
  \omega^{(r)}\cdot \vet p= \tilde{\vet \omega}\cdot \vet p,
\end{equation}
where $\frac 1 2 C^{(r)}\vet p \cdot \vet p$ is the part in normal
form that is quadratic in $\vet p$, i.e., $\sum_{s =0}^{r-1} \langle
f_2^{(r-1,s)}\rangle_{\vet q}$. Moreover, the $r$-th normalization
step is given by $\hat{\Kscr}^{(r)} = \exp
\lie{\chi_1^{(r)}}\Kscr^{(r-1)}$ and $\Kscr^{(r)} = \exp
\lie{\chi_2^{(r)}}\hat{\Kscr}^{(r)}$, where the first generating
function is taking into account also the translation, i.e.,
$\chi_1^{(r)} = X^{(r)}+ \vet \xi^{(r)}\cdot {\vet q}$, while the
solution for $X^{(r)}$ and $\chi_2^{(r)}$ is still given by
equations~\eqref{frm:X}--\eqref{frm:chi2}, where $\vet \omega^{(r-1)}$
is replaced by $\tilde {\vet \omega}$.

Let us suppose to be able to iterate it {\it ad infinitum} the
algorithm, then we would end up with an Hamiltonian of the form
\begin{equation}
\label{eq:haminfty}
\Kscr^{(\infty)}(\boldsymbol{p},\boldsymbol{q})\,=\, %
\tilde{\boldsymbol{\omega}}\cdot\boldsymbol{p} + \Oscr(\|\vet{p}\|^2)\ 
\end{equation}
and, by writing the corresponding equations of motion, it is easy to
realize that the torus $\{\vet{p} = \vet{0}\,,\ \vet{q}\in\toro^2\}$
would be obviously invariant and the motion on it would be
quasi-periodic and characterized by an angular velocity vector equal
to $\tilde {\vet \omega}$. As in the case of the algorithm for the
normal form for elliptic tori, we can explicitly iterate the algorithm
only up to a finite normalization step and we can investigate the
convergence of the procedure by controlling the decrease of the norm
of the generating functions. Let us emphasize that, when the initial
approximation is not good enough, the translations $\vet
\xi^{(r)}\cdot \vet q$ introduced to keep the frequency fixed are not
close to the identity and, therefore, the convergence of the procedure
is prevented. This is the reason why we iterate the algorithm with
fixed frequencies only at the end of the whole application of the
Newton method, when we have approached better the aimed solution.

\subsection{Computer-assisted proof of existence of KAM tori}
\label{subsec:CAP}
We have found interesting to enqueue to the construction of the normal
form a scheme of estimates in order to produce a computer-assisted
proof of the existence of the specific invariant torus we are aiming
at. For what concerns the KAM theorem, let us recall that there are
different computer--assisted proofs in literature (see,
e.g.,~\cite{Cel-Chi-2007} and~\cite{Loc-Gio-2000}), that have been
introduced introduced with the purpose of filling the gap between
analytical theory and applications to Celestial Mechanics. Here, we
follow the approach very recently described in~\cite{Val-Loc-2021} and
in subsection~\ref{sbs:rigorous-results} we will also show the results
produced by using the codes included into the supplementary material
related to that work. However, instead of simply claiming that there
is a software package which can be used as a black-box in order to
produce a computer--assisted proof of existence of KAM tori, we think
that it is more convenient to briefly recall that approach and discuss
its underlying ideas. This is made in our framework also with the
purpose of clarifying how to implement such a computer--assisted proof
with very concrete chances of success in other contexts.

In order to apply an appropriate version of the KAM theorem, let us
define the Hamiltonian for which we want to prove the existence of the
invariant torus by iterating one last time the Kolmogorov
normalization algorithm with a proper scheme of estimates. Our generic
$r-1$-th Hamiltonian now writes
\begin{equation}
  \label{frm:ham-barr-kamtor}
  \Kscr^{(r-1)} = E^{(r-1)} +\tilde{ \vet \omega}\cdot \vet p + \sum_{s\ge
    0} \sum_{\ell \ge 2}f_\ell^{( r-1, s)}(\vet p, \vet q) + \sum_{s\ge
    r} \sum_{\ell = 0}^1 f_\ell^{(r-1, s)}(\vet p, \vet q)\ ,
\end{equation}
since here we aim to iterate the algorithm for constructing the
Kolmogorov normal form in the version keeping fixed the angular
velocity vector $\tilde{ \vet \omega}$. The remainder $\Rscr^{( r-1)}
= \sum_{s\ge r} \sum_{\ell = 0}^1 f_\ell^{(r-1, s)}(\vet p, \vet q)$
is the quantity we are interested in estimating and reducing below the
threshold value of applicability.  Following the approach used
in~\cite{Car-Loc-2020} and~\cite{Val-Loc-2021} which goes back to the
technique developed in~\cite{Cel-Gio-Loc-2000}, we distinguish between
terms for which we will explicitly compute the Kolmogorov normal form
(with the help of an algebraic manipulator,
see~\cite{Gio-San-Chronos-2012} for an introduction) up to $R_{\rm
  I}$-th step, terms for which we estimate the
$\ell_1$-norms\footnote{Let us recall that in the present framework
  the $\ell_1$-norm for real functions having a finite representation
  in Taylor-Fourier series is given by $\|f_l^{(r-1,r)}\| =
  \sum_{j_1+j_2=l} \sum_{0\le|\vet k|\le 2s} \big(|c_{\vet j, \vet
    k}^{(r-1,s)}|+|d_{\vet j, \vet k}^{(r-1,s)}|\big)$ being
  $f_l^{(r-1,r)}(\vet p, \vet q) = \sum_{j_1+j_2=l} \sum_{0\le|\vet
    k|\le 2s} p_1^{j_1}p_2^{j_2}\big[c_{\vet j, \vet
      k}^{(r-1,s)}\cos({\vet k \cdot \vet q})+d_{\vet j, \vet
      k}^{(r-1,s)}\sin({\vet k \cdot \vet q})\big]$ the corresponding
  expansion of the function.\label{nota:def-norma}} with constants
$\Fscr_\ell^{(r,s)}$ such that $\|f_\ell^{(r,s)}\|\le
\Fscr_\ell^{(r,s)}$ (up to a maximal degree in the angles $\vet q$
that is equal to $2R_{\rm II}$) and the tail of infinite terms. In the
case of the latter ones, we provide just an uniform estimate
\begin{equation}
\label{frm:stimacode}
\|f_\ell^{(r-1,s)}\|\le a_{r-1}^{s}\, \zeta_{r-1}^\ell\, \Escr_{r-1}\ ,
\end{equation}
being $a_{r-1}$, $\Escr_{r-1}$ and $\zeta_{r-1}$ suitable constants.
This means that at every step $r-1$ the Hamiltonian can be completely
represented by the following finite list of quantities:
\begin{equation}
\label{frm:hamCAP}
\vcenter{\openup1\jot\halign{
 \hbox {\hfil $\displaystyle {#}$}
&\hbox {\hfil $\displaystyle {#}$\hfil}
&\hbox {$\displaystyle {#}$\hfil}\cr
 \Sscr^{(r-1)} &= \bigg\{
 & \tilde{\vet \omega},\ f_{0}^{(r-1,\min\{r, R_{\rm I}\})} \,,\,\ldots\,,f_{0}^{(r-1,R_{\rm I})},  f_{1}^{(r-1,\min\{r, R_{\rm I}\})} \,,\,\ldots\,,f_{1}^{(r-1,R_{\rm I})},\cr
 & & f_{2}^{(r-1,0)} \,,\,\ldots\,,f_{2}^{(r-1, R_{\rm I})},
 \cr
 & &\Fscr_{0}^{(r-1, R_{\rm I}+1)} \,,\,\ldots\,, \Fscr_{0}^{(r-1, R_{\rm II})}  \,,\, 
\Fscr_{1}^{(r-1,R_{\rm I}+1)} \,,\,\ldots\,, \Fscr_{1}^{(r-1, R_{\rm II})}  \,,\,\cr
& &\Fscr_{2}^{(r-1, R_{\rm I}+1)} \,,\,\ldots\,, \Fscr_{2}^{(r-1, R_{\rm II})}  \,,\,
 \cr
 & & \Escr_{r-1} \,,\,  a_{r-1},\, \zeta_{r-1}  \bigg\}
 \ .
 \cr
}}
\end{equation}
where, for the sake of brevity, we have reported just the Hamiltonian
terms up to degree $2$ in the actions, even if it can be extended to a
general finite degree. Let us emphasize that all the coefficients that
appear in the Taylor-Fourier expansions of the functions
$f_\ell^{(r-1,s)}$ (making part of the previously described set
$\Sscr^{(r-1)}$) are computed by using interval arithmetic. Moreover,
in order to avoid a too fast growth of the constants
$\Fscr_\ell^{(r,s)}$ and the exceeding of the limits in the size of
the numbers which can be safely represented on a computer, we do the
practical computation of these values by using their logarithms.  The
main difficulty is therefore to understand how to redefine the
quantities in $\Sscr^{(r)}$, i.e., after an application of the change
of coordinates for the construction of the normal form as it has been
described at the end of the previous subsection~\ref{sec:kamtor}.  The
essential effort is required by the evaluation of the effects induced
by the Lie series. More details about the technical estimates needed
for the proof are reported in appendix~\ref{sec:technicalities}. The
use of validated numerics (see the appendixes of~\cite{Car-Loc-2020}
for a gentle introduction to some of its main concepts) in all the
prescribed computations can make the computer-assisted proof
completely rigorous.

The final part of the computer-assisted proof is devoted to the check
that all the hypotheses of the KAM theorem in the statement
reported\footnote{Such a particular version of the KAM statement has
  been chosen, because in~\cite{Stef-Loc-2012} the threshold of
  applicability of the theorem can be explicitly evaluated in a way
  fitting well with our computer-assisted scheme of estimates.}
in~\cite{Stef-Loc-2012} are satisfied for a Hamiltonian
$\Kscr^{(R_{\rm II})}$, where $R_{\rm II}$ is the final number of the
normalization steps for which we can explicitly provide the list of
quantities appearing in the set $\Sscr^{(R_{\rm II})}$. This requires
to suitably translate the upper bounds on the $\ell_1$-norm in
estimates on the sup-norm of analytic functions defined on a domain
$D_\rho\,$, being $D_\rho = \{(\vet p, \vet q) \in \complessi^{2n} :
\max_j| p_j | <\rho,\ \max_j | {\rm Im }(q_j) |<\rho\}$. Of course, we
restrict to consider the set of all the analytic functions defined on
$D_\rho$, that have period $2 \pi$ in $q_1,\ldots,q_n$ and are real
for real values of the variables. Once the quantities listed in the
set $\Sscr^{(R_{\rm II})}$ are explicitly known for $R_{\rm II}$ large
enough, providing the needed upper bounds on the sup-norms does not
require a lot of additional technical work as it is explained in
detail at the end of section~3.3 of~\cite{Val-Loc-2021}. This is true
also for both the check of the non-degeneracy condition on the Hessian
of the quadratic part of the Hamiltonian (so that an equation similar
to that in~\eqref{frm:trasl} can {\it always} be solved) and the
complete determination of a pair $(\tilde{\vet \omega},\gamma)$, being
$\tilde{\vet \omega}$ a vector that satisfies a Diophantine inequality
of type~\eqref{diseq:Dioph}, where the constant value of parameter
$\gamma$ appears.

\section{Results based on a normal forms constructive approach}
\label{sec:results}
In this section we will discuss in a more quantitative way the results
obtained by following the procedures described in
Section~\ref{sec:toro} for the planetary system
$\upsilon$~Andromed{\ae}. More precisely, we refer to its three-body
planetary model as it is completely defined by the initial conditions
and the parameters listed in Table~\ref{tab:dopo-griglia}.

\subsection{Secular dynamics: comparisons between numerical integrations and semi-analytic computations}
\label{sbs:semianalytical-results}

First, we focus on the construction of a suitable lower dimensional
elliptic torus, that is a milestone on the way to reach the wanted
librational KAM torus, as it has been widely discussed in
Section~\ref{sec:toro}. For what concerns the initial secular model
defined in formula~\eqref{eq:hamsecazang}, we truncate its expansion
at degree~$8$ in $(\sqrt{J_1},\sqrt{J_2})$, while all the Hamiltonians
of type $H^{(r)}$ that are introduced in section~\ref{subsec:ell} have
been expanded up to degree 8 in the actions. During the iteration of
the normalization steps constructing the lower dimensional elliptic
torus, our computational codes are representing Hamiltonian terms with
a trigonometrical degree up to $72$ in the angles. This choice is
consistent with the number of steps performed for both the
normalization algorithms (the one for the construction of the elliptic
torus and the one for the KAM torus), that is $35$. Let us recall that
at the $r$-th step the generating functions $\chi_0^{(r)}$,
$\chi_1^{(r)}$, $X_2^{(r)}$ and $Y_2^{(r)}$ are determined in such a
way to remove the perturbing terms of trigonometrical degree up to
$2r$.

\begin{figure}
\centering
\includegraphics[width=7.9cm]{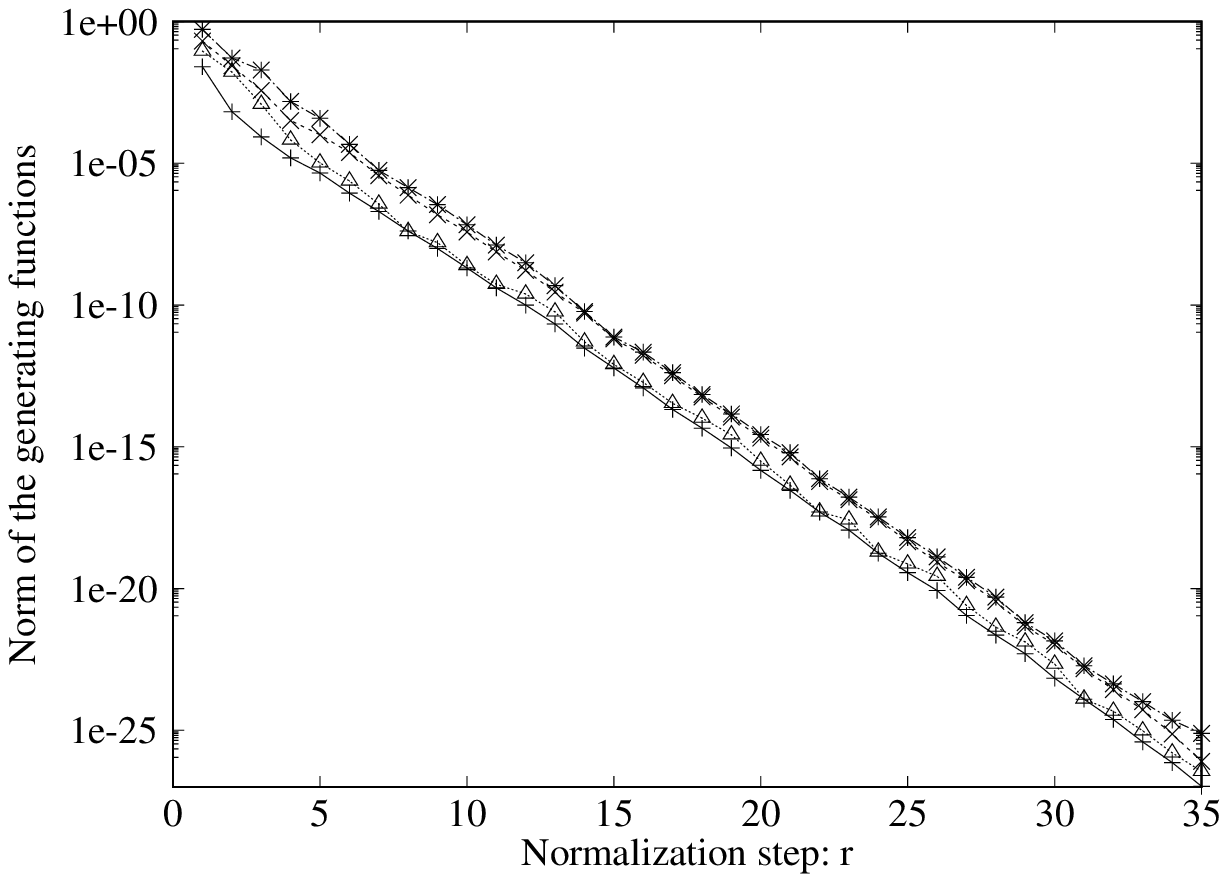}
\includegraphics[width=7.9cm]{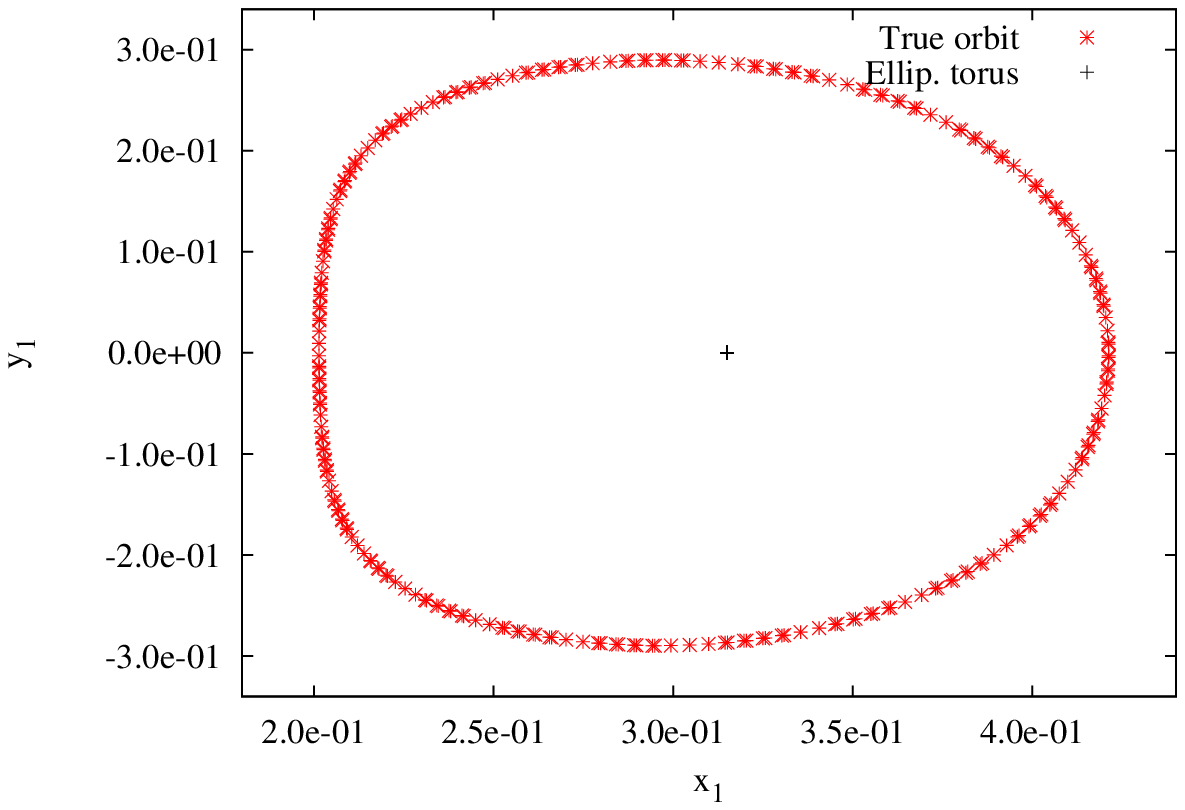}
\caption{On the left, norms (in log scale) of the generating functions
  defined in the algorithm for the construction of elliptic tori. On
  the right, the Poincar\'e sections of the constructed periodic
  solution, enclosed in the {\it true secular orbit} that is plotted
  as well.}
\label{fig:normegen-ell}
\end{figure}

In the left panel of Figure~\ref{fig:normegen-ell}, we plot the norms
of the generating functions (with the exception of the one responsible
of the diagonalization) introduced by the algorithm for the
construction of the normal form for elliptic tori. Said norm is
defined as the sum of the absolute values of every coefficient
appearing in the expansions.  All the norms decrease regularly and
very quickly. Let us emphasize that the same algorithm
in~\cite{Car-Loc-2021} has given the best results in the FPU
$\beta$-model, where, like in this case, the perturbation was quartic,
with an Hamiltonian that is an even polynomial in the canonical
variables. Moreover, being this model a purely secular approximation
of a planetary system, there is not any distinction between fast and
slow dynamics, as in the framework considered in~\cite{Car-Loc-2021}
and differently from the planetary problem considered
in~\cite{Gio-Loc-San-2014}.  At the end of the normalization we obtain
an Hamiltonian as in~\eqref{frm:ham-barr-elltor}. Therefore,
neglecting the remainder we can define a set of initial conditions for
a periodic orbit by taking $(p(0),q(0),x(0),y(0)) = (0,0,0,0)$. We
expect that the motion should be given by $(p(t),q(t),x(t),y(t))
\simeq (0, \omega^{(\bar r)}t,0,0)$. As a test of our result, we can
express the initial conditions in any set of canonical coordinates
that is preliminary to those introduced by the computational algorithm
constructing the normal form~\eqref{frm:ham-barr-elltor}, which
corresponds to (our approximation of) an elliptic torus.  For
instance, this can be done for the canonical variables $(\vet x, \vet
y)$ that are defined by the equations~\eqref{eq:azang}
and~\eqref{eq:def-azang-Iphi}--\eqref{eq:def-dalembcoord-xy} starting
from the initial coordinates $(\vet \xi, \vet \eta)$ entering into the
definition of the purely secular model~\eqref{frm:hsec}.  The
procedure computing unknown values of canonical variables by
composition of several canonical transformations has been thoroughly
described in several works using the Lie series technique, in
particular we can defer the interested reader to section 4.2
of~\cite{San-Loc-Gio-2011}. We can therefore numerically integrate the
equations of motion related to the Hamiltonian model~\eqref{frm:hsec},
starting from $(\vet \xi(0), \vet \eta(0))$ which corresponds to
$(p(0),q(0),x(0),y(0)) = (0,0,0,0)$; we can then compute its
Poincar\'e sections with respect to the $(x_1\,,\,y_1)$--plane.
Looking at the right panel of Figure~\ref{fig:normegen-ell}, we can
appreciate that the orbit in black appears as a fixed point and that,
therefore, the periodic solution is well approximated.  The accuracy
of the result can be checked also by comparing the frequencies of the
numerical elliptic torus with the one constructed with our algorithm:
the discrepancy is of order $10^{-7}$, i.e., the same order of the
tolerated error in the identification of the elliptic torus through
the FA.  Let us emphasize that, considered the fast decrease of the
perturbation, $35$ steps are way more than necessary in order to
produce a valuable approximation of the elliptic torus, which is not
the final aim of this work. Nevertheless, we take this value because
of the relatively small amount of time needed for the
computation\footnote{About one hour of CPU-time was requested in order
  to perform such a computation on a modern workstation equipped with
  a \texttt{Xeon 18-Core 6154 - 3.0Ghz}.\label{nota:CPU-specs}} and in
view of the fact that we are going to use the same truncations for
what concerns the following Kolmogorov algorithm, where we will see
that the convergence of the procedure is not equally fast.  We recall
also that the iteration of the algorithm introduces corrections on the
value of the energy and on the frequencies of the periodic orbit, so
that the iterations of the Newton methods are essential in order to
converge to the same periodic orbit numerically found.

The composition of canonical transformations can be used to express
the initial conditions of the {\it true orbit} as a function of the
canonical variables introduced to define the normal form of the
elliptic torus.  In particular, we can compute the initial value of
the action $I_2^\star= (x^{\star2} + y^{\star2})/2$ at time $t=0$ for
the {\it true orbit} in the variables of the normal form. Let us
underline that the value of the translation $I_2^\star$ is defined as
the area delimited by a section of the orbit, orthogonal to the
elliptic torus. The smaller the translation, the closer the {\it true
  orbit} should be to the elliptic torus and, therefore, in a less
perturbed zone. Here, we found convenient to choose the periodic orbit
whose period is close to the one of the {\it true orbit}, as that is
the first approximation used in the iteration of the FA for locating
the lower dimensional torus.

\begin{figure}
\centering
\includegraphics[width=7.9cm]{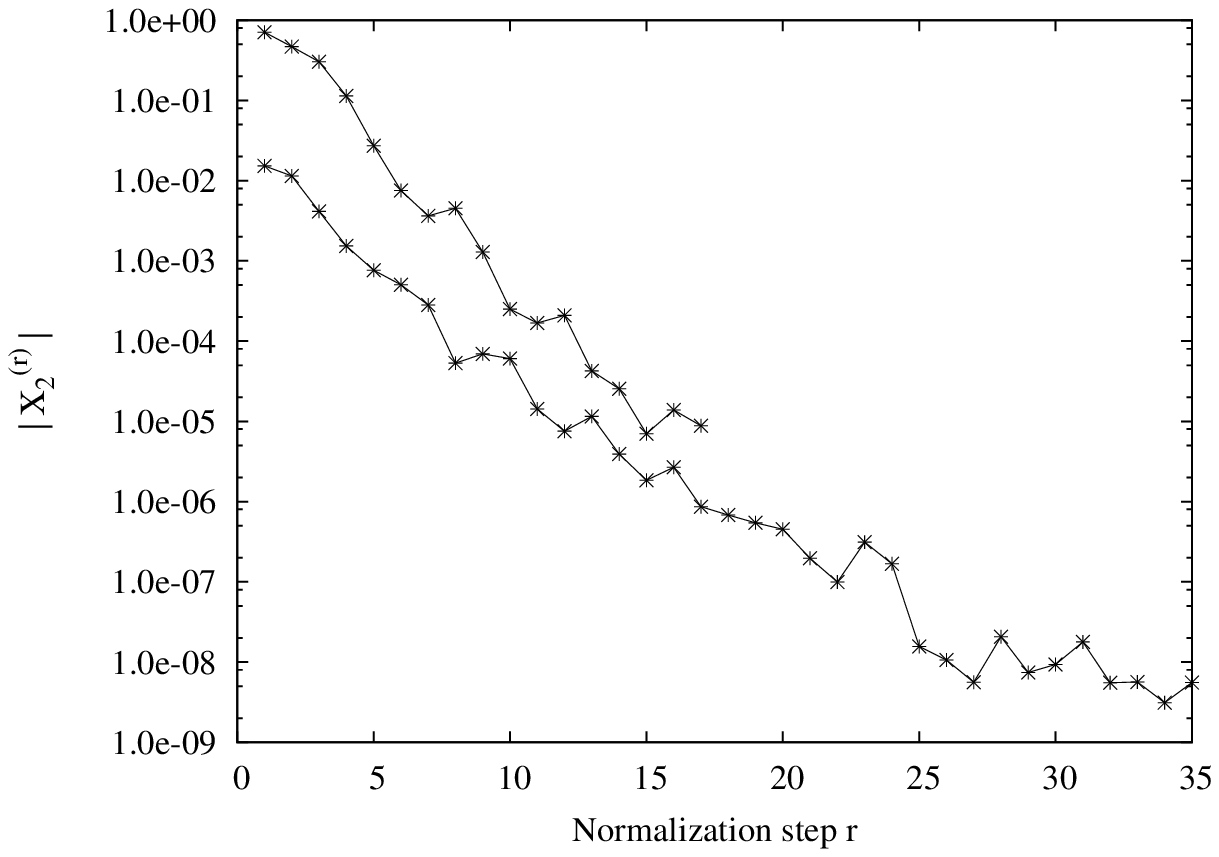}
\includegraphics[width=7.9cm]{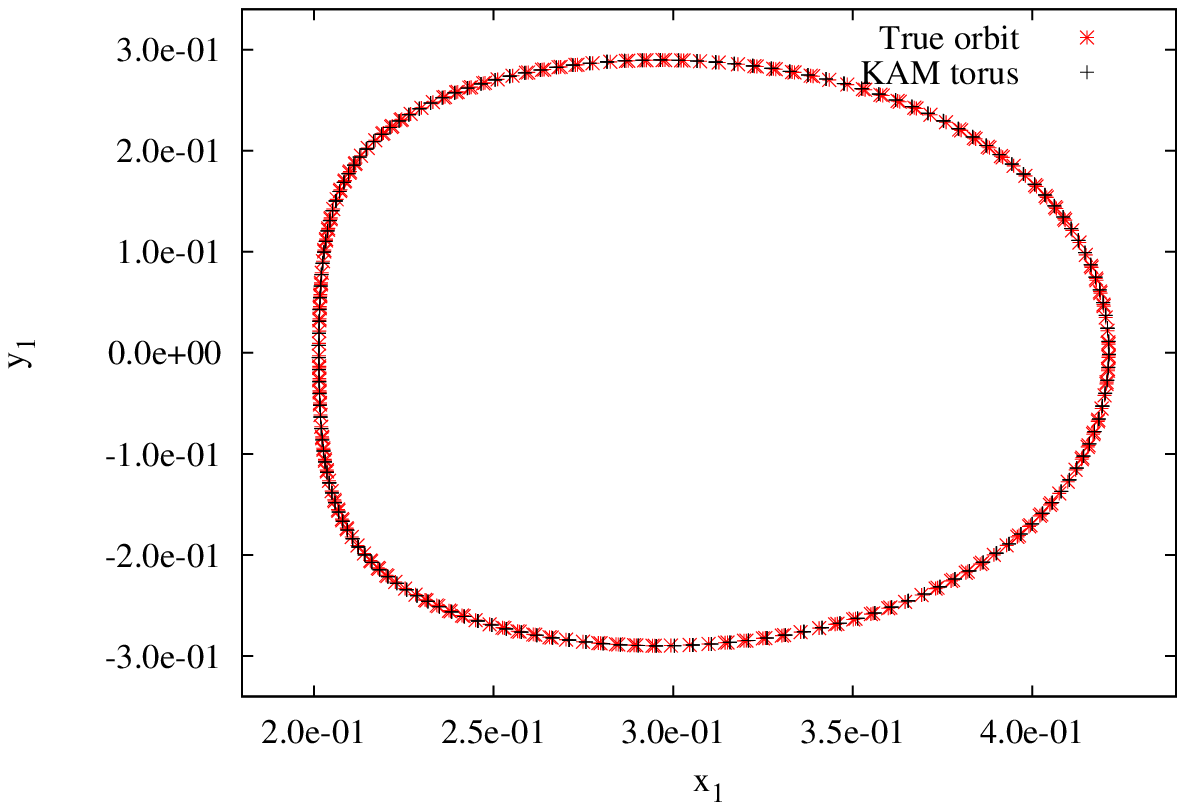}
\caption{On the left, the norms (in log scale) of the generating
  functions $\chi_2^{(r)}$ introduced by the algorithms for the
  construction of the normal form for a KAM torus: the one above
  refers to the algorithm without fixing the frequencies and the one
  below to the algorithm with the translation for keeping the
  frequency vector fixed as $\tilde{\vet \omega}$. On the right, the
  Poincar\'e sections of the {\it true orbit} (in red) and of the
  orbit with initial conditions on the KAM torus constructed with our
  algorithm (in black).}
\label{fig:normegen-kam}
\end{figure}

At this level, we start to iterate the algorithm for the construction
of the normal form for a KAM torus, in a first version that is without
introducing any other translation to target the frequencies of the
quasi-periodic motion characterizing the wanted invariant
manifold. This is equivalent to construct one of the intermediate tori
between the elliptic torus and the {\it true orbit}. By using the
Newton method described at the end of the previous
Section~\ref{sec:toro}, this computational procedure can converge to
values $(\tilde I_1\,,\,\tilde I_2)$, that enter into the
definition~\eqref{solo1trasl} of the initial translation $\chi(\vet
q)$. These values are such that the relative error between the
frequency of the approximated KAM torus $\vet \omega^{(r)}$ and the
numerical one $\tilde{\vet \omega}$ is less than $10^{-3}$. Then, we
apply one last time the classical KAM algorithm, i.e., with the
addition of a translation at each normalization step in order to keep
the frequency fixed. In Figure~\ref{fig:normegen-kam}, we plot the
norms of the generating functions $\chi_2^{(r)}$ for both the
normalization algorithms, the one without any other translation apart
from the initial one $\chi(\vet q)$ defined in~\eqref{solo1trasl} and
the one enqueued with the translation $\vet \xi^{(r)}\cdot \vet q$ at
each step $r$ for aiming at the wanted frequency vector
$\tilde{\vet\omega}$. We can observe in both cases a slow decrease of
the norm and, in particular, that we have an improvement in iterating
the algorithm at fixed frequency only subsequently, when the
perturbation has been considerably reduced. Being the convergence rate
quite slow, it is not so rare to encounter high order resonances that,
even if they do not compromise the convergence, have an evident impact
in the norms of the generating functions, appearing as small and
sudden peaks. This is the reason why we have iterated the first
algorithm just half the steps in order to reduce the danger of falling
in a resonant region because of the random variation of the sequence
$\vet \omega^{(r)}$. Such a problem is not truly related to the
non-resonance properties of the frequency vector $\tilde{\vet\omega}$.

In the same way as before, we can evaluate the composition of all the
transformations of coordinates introduced into the framework of the
previous Section~\ref{sec:toro}, with the aim of constructing a well
defined sequence of Hamiltonians, that in principle should converge to
the Kolmogorov normal form~\eqref{eq:haminfty}. Therefore, we can
consider the initial conditions $(\vet \xi(0), \vet \eta(0))$ which
corresponds to a point $({\vet p}(0),{\vet q}(0)) = ({\vet 0},{\vet
  0})$ belonging to the torus that we have built and is expected to be
invariant with a very good approximation.  Hence, we integrate
numerically the equation of motions starting from the initial
conditions $(\vet \xi(0), \vet \eta(0))$ and we can plot the
corresponding Poincar\'e sections with respect to the
$(x_1\,,\,y_1)$--plane. In the right panel of
Figure~\ref{fig:normegen-kam}, we compare the Poincar\'e sections of
this orbit with those of the {\it true secular orbit}, obtaining an
excellent superposition. Such a comparison confirms that the {\it true
  orbit} should lie very close to a KAM torus, because it can be
approximated very well by the solution we have constructed using a
normal form technique.

\subsection{Rigorous results}
\label{sbs:rigorous-results}

By comparing the two plots in the left panels of
Figures~\ref{fig:normegen-ell}--\ref{fig:normegen-kam}, one can
appreciate that the decrease of the norms in the final Kolmogorov
algorithm is way slower and more irregular with respect to the case of
the elliptic torus. This remark makes doubtful the convergence to a
Kolmogorov normal form for what concerns the sequence of Hamiltonians
introduced in subsection~\ref{sec:kamtor}. In this framework, the most
convincing argument is definitely provided by a rigorous
computer-assisted proof, that we will discuss in the following.

According to the strategy designed in~\cite{Cel-Gio-Loc-2000}, it is
convenient to preliminarly perform a large enough number of
normalization steps constructing the Kolmogorov normal form. This is
made with the purpose of reducing the size of the perturbative terms
so much that it is possible to apply a suitable version of the KAM
theorem. In practice, the perturbation is made smaller by applying a
procedure where a couple of different stages can be distinguished.
First, the expansions of the generating functions $X^{(r)}$,
$\vet\xi^{(r)}\cdot \vet q$ and $\chi_2^{(r)}$ are explicitly computed
when the index $r$ of the normalization step is running from $1$ to
$R_{\rm I}\,$; then, just the estimates of the norms of terms up to
order $\epsilon^{R_{\rm II}}$ are iterated, by avoiding the
calculation of their cumbersome expansions if they are
$\Oscr(\epsilon^{r})$ with $R_{\rm I}<r\le R_{\rm II}\,$. All the
terms appearing in a Hamiltonian $\Kscr^{(r-1)}$ satisfy suitable
inequalities that involve a finite number of quantities in the
definition of the upper bounds, being the latter ones collected in a
set $\Sscr^{(r-1)}$ of the type described
in~\eqref{frm:hamCAP}. Appendix~\ref{sec:technicalities} includes the
recursive formul{\ae} allowing to properly define the values of the
elements that belong to the set~$\Sscr^{(r)}$ and provide the upper
bounds for the terms making part of Hamiltonian $\Kscr^{(r)}$. In
particular, the way to iterate the estimates of the norms of the
generating functions is fully explained in such an appendix.  The
number $R_{\rm I}$ of the first normalization steps that are performed
by explicitly computing the expansion of the generating functions must
be large enough, in order to subsequently iterate the estimates in a
so successful way to reduce the size of the perturbation as much as
needed.

\begin{figure}
\centering
\includegraphics[width=7.9cm]{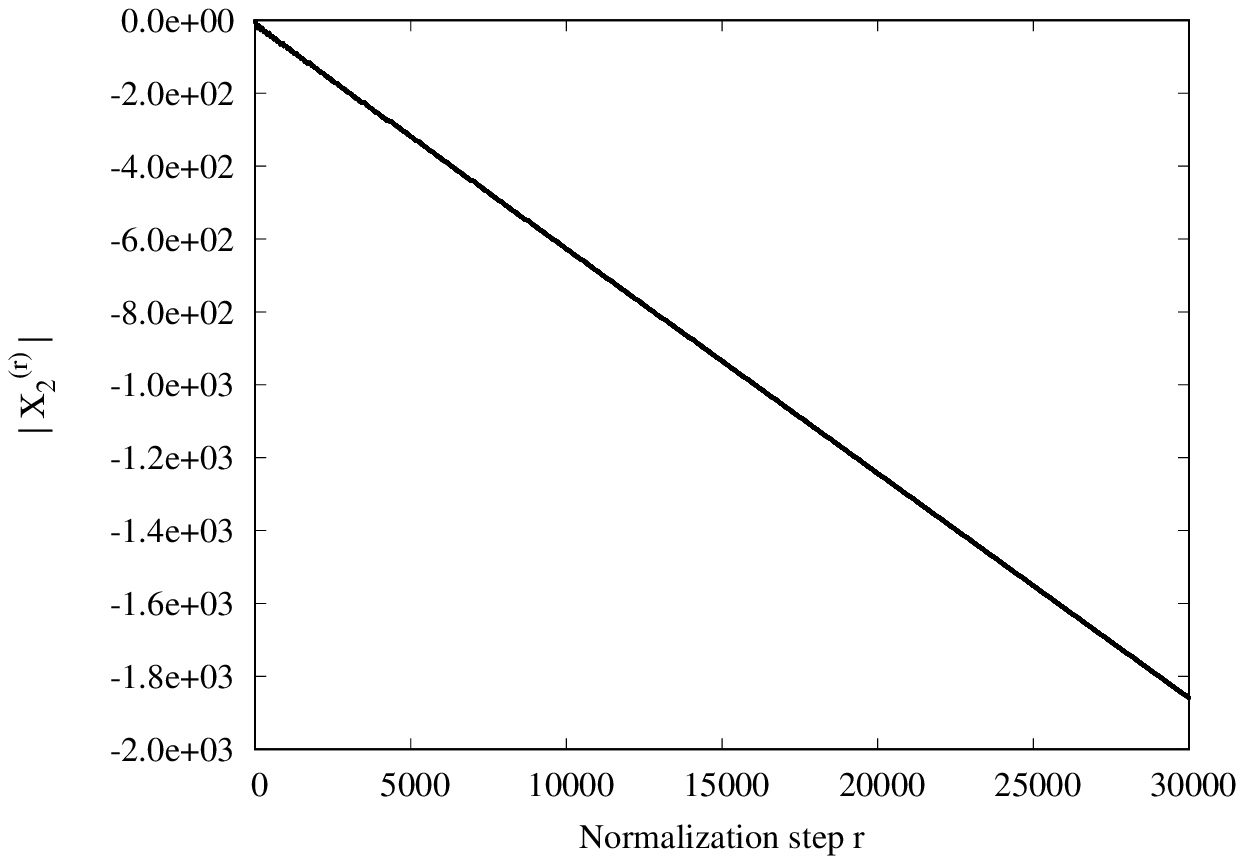}
\includegraphics[width=7.9cm]{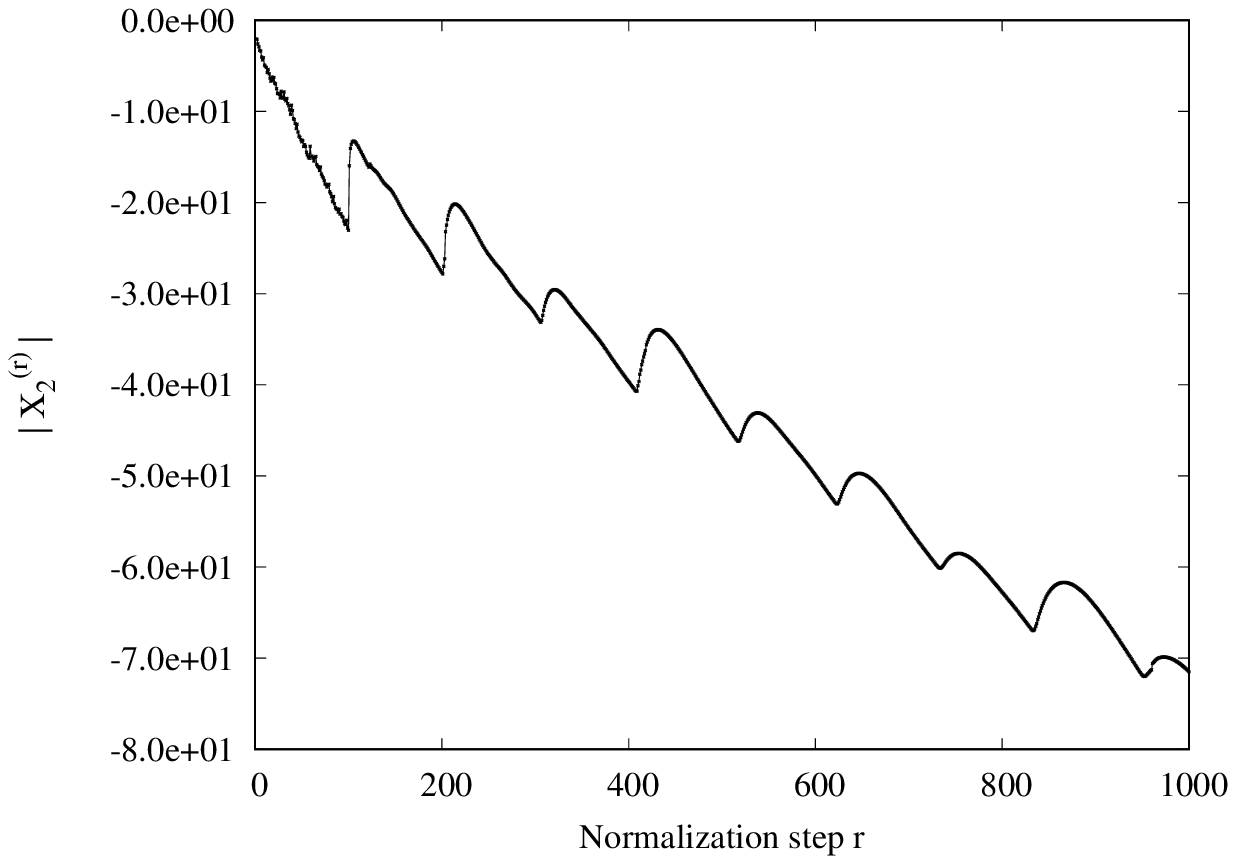}
\caption{Estimates of the norms (in semi-log scale) of the generating
  function $\chi_2$. On the right, the zoom on the first 500
  normalization steps.}
\label{fig:stime}
\end{figure}

For our purposes, we have seen that it is enough to set $R_{\rm I}=60$
and to represent the estimates of the norms of terms up to order
$\epsilon^{R_{\rm II}}$ with $R_{\rm II} = 10000$. The plots of the
estimates of the norm of $\chi_2^{(r)}$ are reported in
Figure~\ref{fig:stime} and confirm the slow but regular reduction of
the perturbation. Indeed, the plot in the right box of such a figure
highlights that the decrease is steeper up to the step $60$, i.e.,
when an explicit computation of the generating function is performed;
then, we have a sudden leap, due to the transition to the regime of
the pure iteration of the estimates; thereafter, a periodic and
increasingly small jump is noticeable until the decrease stabilizes to
a flatter but nearly constant slope.  As it is discussed also in the
explanatory document included in the software package allowing to
implement such a kind of computer-assisted proofs
(see~\cite{Loc-CAP4KAM-2021}), it is an important achievement when the
behavior of the norms of the generating functions looks similar to
that depicted in Figure~\ref{fig:stime}. In fact, this ensures that
the existence of quasi-periodic motions of frequency vector
$\tilde{\vet \omega}$ and invariant with respect to the Kolmogorov
normal form~\eqref{eq:haminfty} can be rigorously proved, by
eventually increasing the index of the final normalization step
$R_{\rm II}\,$. In the case of the model we are studying, the
perturbative terms appearing in the expansion of the Hamiltonian
$\Kscr^{(10000)}$ are so small that the computer-assisted proof can be
successfully completed by applying the KAM theorem in the version
described in~\cite{Stef-Loc-2012}.  Performing the first $R_{\rm
  I}=60$ normalization steps so as to explicitly determine the expansion
of all the corresponding generating functions requires nearly 7 hours
on the same computer whose CPU is described in
footnote${\ref{nota:CPU-specs}}\atop{\phantom{1}}$; moreover, about
two additional hours of CPU-time are needed by the iteration of the
estimates for the norms of terms up to order $\epsilon^{R_{\rm II}}$
with $R_{\rm II}=10000$. Of course, the larger amount of time that is
requested by the explicit construction of the normal form can be
easily explained because of the huge number of terms to
manipulate. This is also due to our decision to represent in the
expansions of the finite sequence of Hamiltonians
$\big\{\Kscr^{(r)}\big\}_{r=0}^{60}$ all the terms having up to degree
$8$ in the actions, while the maximal trigonometrical degree has been
fixed to $2R_{\rm I}=120$.  Such a truncation rule on the actions is
in agreement with the limits we have adopted during the the
computation of the Hamiltonians $\Hscr^{(0)},\ldots,\,\Hscr^{(17)}$,
that has been made according to the algorithm described in
subsection~\ref{sec:kamtor}.  These preliminary expansions have been
basic to produce the semi-analytical results already discussed in the
previous subsection~\ref{sbs:semianalytical-results}.

The result we have obtained thanks to a rigorous proof can be summarized as
follows.
\begin{theorem}[Computer-assisted]
Let us consider the Hamiltonian $\Kscr^{(0)}=\Hscr^{(17)}$ whose
expansion is written in~\eqref{frm:hbarr} that has been obtained by
following the procedure described in
sections~\ref{subsec:ell}-\ref{sec:kamtor} and by truncating all the
expansions up to degree $8$ in the actions and to the trigonometrical
degree $72$ in the angles.  Let $\tilde{\vet\omega}\in\reali^2$ be
such that
\begin{equation}
\vcenter{\openup1\jot\halign{
 \hbox {\hfil $\displaystyle {#}$}
&\hbox {\hfil $\displaystyle {#}$\hfil}
&\hbox {$\displaystyle {#}$\hfil}\cr
 \tilde{\omega}_1 &\in
 &( -1.3375038001502735\times 10^{-2}\,, \ -1.3375038001500735 \times 10^{-2})\,,
 \cr
 \tilde{\omega}_2 &\in
 &(\ 7.4969769908765335 \times 10^{-3}  \,,\ 7.4969769908785335\times 10^{-3})\,,
 \cr
}}
\label{rect-bounds:omega-CAP}
\end{equation}
and it satisfies the Diophantine\footnote{Both the midpoints and the
  widths of the intervals reported in
  formula~\eqref{rect-bounds:omega-CAP} are compatible with the values
  of the angular velocities that we have obtained by using the methods
  of the frequency analysis (see, e.g.,~\cite{Laskar-03}). Such a
  technique has been applied to study the numerical integration of the
  equations of motion corresponding to the Hamiltonian
  model~\eqref{frm:hsec}, starting from the initial conditions
  corresponding to the orbital elements reported in
  Table~\ref{tab:dopo-griglia}. At the end of the computer-assisted
  proof, the code that is in charge to check for the applicability of
  the KAM theorem also ensures that in the Cartesian product of those
  intervals, there exist Diophantine vectors satisfying
  inequality~\eqref{diseq:Dioph}. In particular, this is done for the
  very peculiar subset of the pairs
  $\big(\tilde{\omega}_1\,,\,\tilde{\omega}_2\big)$ that are bounded
  so as to belong to the rectangular box~\eqref{rect-bounds:omega-CAP}
  and their ratio is the following algebraic number
  $\big[1523989\,(\sqrt{5}+1)/2\,+\,305410\big] \,/\,
  \big[2718884\,(\sqrt{5}+1)/2\,+\,544869\big]\,$. We have decided to
  focus on this special kind of 2D-vectors, that are selected by using
  the so called computational algorithm of the ``Farey tree'', because
  they are expected to correspond to KAM tori that are particularly
  robust with respect to the perturbations (see~\cite{Kim-Ost-1986}
  and~\cite{Mac-Sta-92}).}  inequality
\begin{equation}
  |\vet k \cdot \tilde{\vet \omega}| \ge \frac{\gamma}{|\vet k|} \qquad
  \forall \ \vet k \in \interi^2\setminus\{\vet 0\}
  \label{diseq:Dioph}
\end{equation}
with $\gamma = 2.094439 \times 10^{-4}$.

Therefore, there exists an analytic canonical transformation
conjugating $\Hscr^{(17)}$ to the Kolmogorov normal form
$\Kscr^{(\infty)}$, which is written in~\eqref{eq:haminfty}. It is such
that the torus $\{\vet p = \vet 0, \vet q \in \toro^2\}$ is invariant
and travel led by quasi-periodic motions whose corresponding angular
velocity vector is equal to $\tilde{\vet \omega}$.
\end{theorem}

Let us recall that the dynamical stability of our planetary model is
expected to follow from the proof of existence of an invariant torus
close to the {\it true orbit}.  Indeed, two possible approaches can be
successfully adopted in this framework. By following the strategy
designed in~\cite{Loc-Gio-2000}, one could prove the existence of a
KAM torus that encloses the {\it true orbit} in the Poincar\'e
sections represented in the first panel of
Figure~\ref{fig:poin-sec}. Since the Hamiltonian secular model has two
degrees of freedom, the constant energy surface is 3D and so 2D
invariant tori that certainly include the initial conditions are able
to ensure a perpetual topological confinement of the {\it true
  orbit}. Otherwise, a more general approach (that is valid also for
Hamiltonian systems having more than two degrees of freedom) can be
applied so to lead to the effective stability concept: in the
neighborhood of the KAM torus, a Birkhoff normal form can be
constructed and complemented with the classical estimates {\it \`a la}
Nekhoroshev (see~\cite{Mor-Gio-1995}). This allows to prove that the
time needed to eventually escape from a small enough region
surrounding the KAM torus can largely exceed the expected life-time of
such a system. Such an approach can be used to show the effective
stability of Hamiltonian planetary models (see~\cite{Gio-Loc-San-2009}
and~\cite{Gio-Loc-San-2017}). Both these approaches would require a
demanding effort for their actual implementations, that we will avoid
because we believe that a complete proof of the stability of the
secular system goes beyond the aims of the present work.

\subsection{Comparisons between the dynamics for the complete planetary model and the secular one}
\label{sbs:confr-modello-completo}
From a more astronomical point of view, it is natural to raise the
following question: how far is reality from the secular Hamiltonian
model we are considering? There are different possible answers to such
a question. In order to give some insights, it is natural to compare
the solutions for the secular model with those related to the complete
planetary system introduced at the very beginning in
formula~\eqref{frm:ham3cpiniz}. Such a comparison can be performed
just for the motion laws that can be observed in both these dynamical
systems, these obviously refer to the secular variables.

\begin{figure}
\centering
\includegraphics[width=7.9cm]{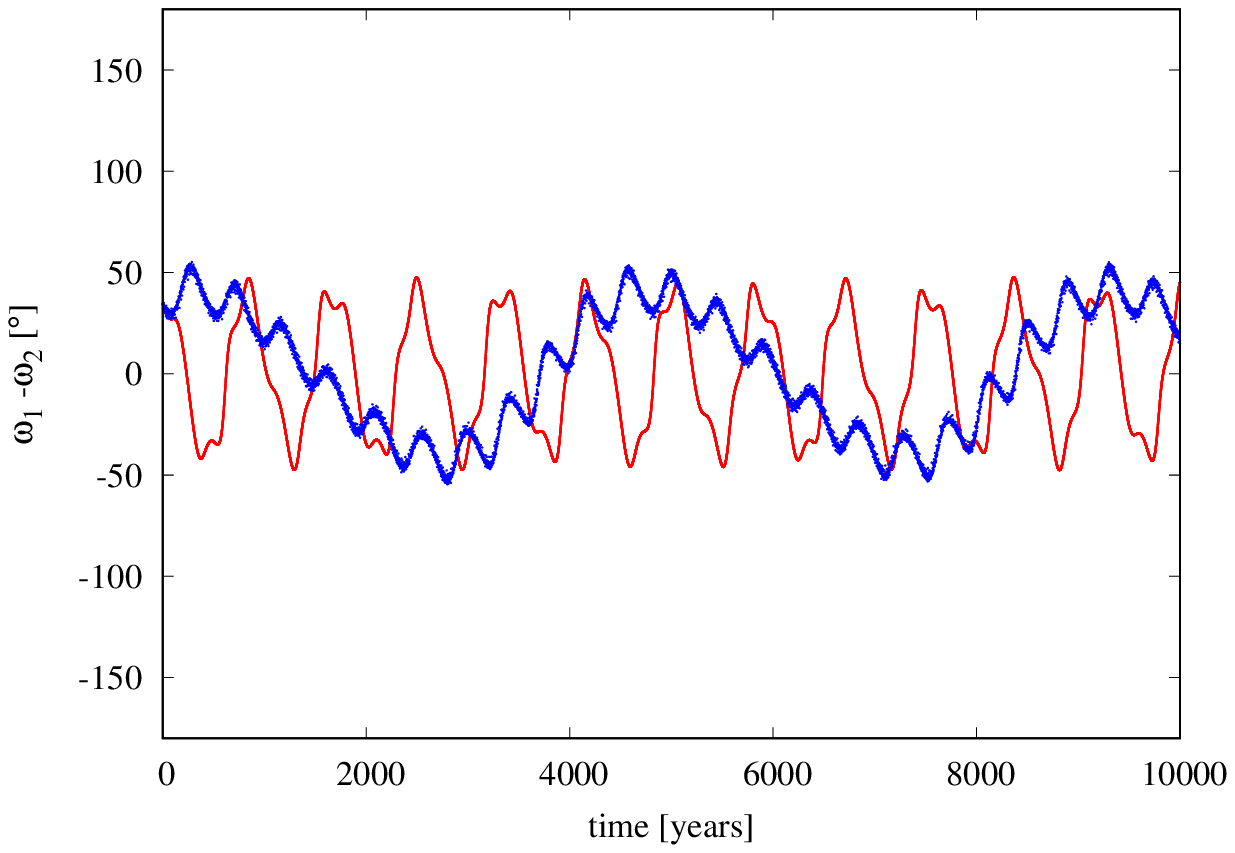}
\includegraphics[width=7.9cm]{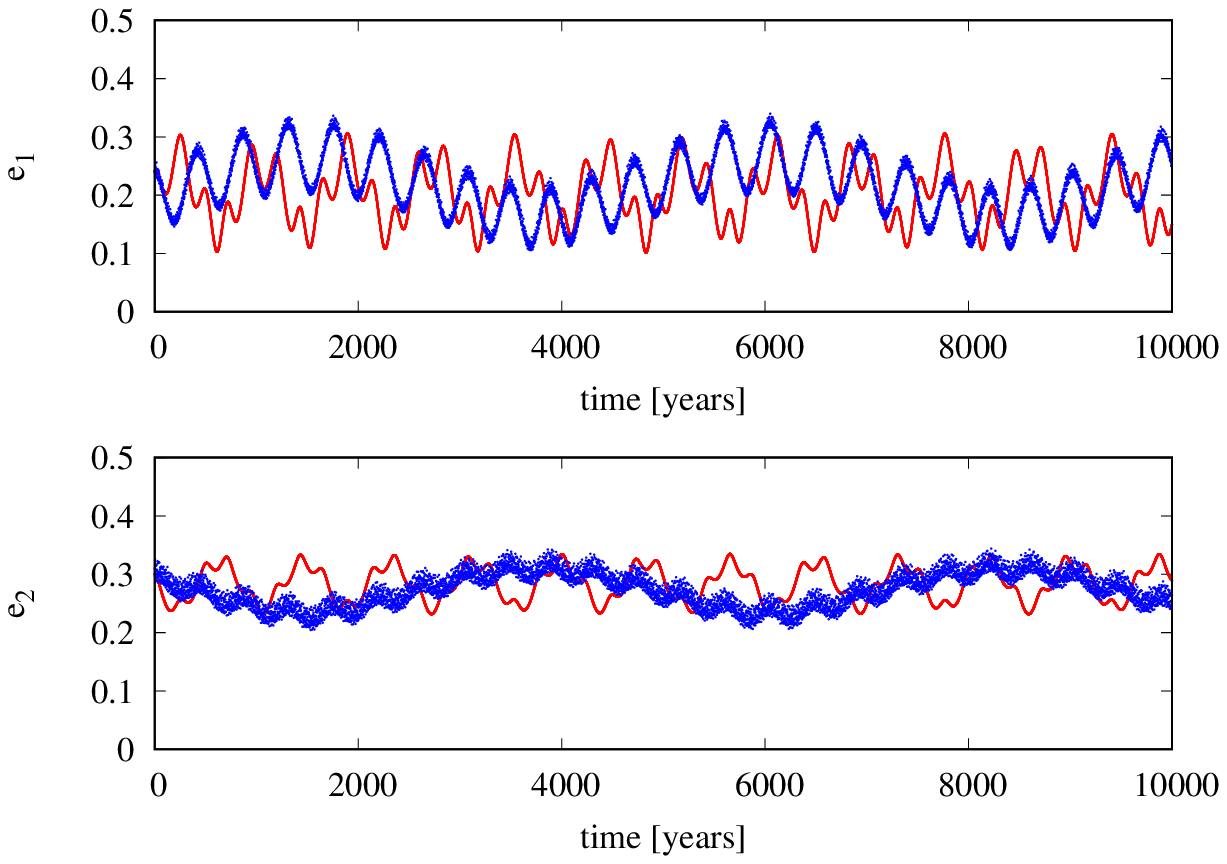}
\caption{Comparisons between the results provided by the numerical
  integrations of the complete planetary system and the secular one
  (in blue and in red, respectively). In the left box, the behavior
  of the difference of the arguments of the pericenters is plotted as
  a function of time. On the right, the temporal evolutions of the
  eccentricities are reported: the plot above [below] refers to
  $\upsilon$~And~\emph{c} [$\upsilon$~And~\emph{d}, resp.].}
\label{fig:secVScom}
\end{figure}

On the one hand, we have produced the numerical integrations of the
equations of motion related to the complete planetary
Hamiltonian~\eqref{frm:ham3cpiniz}, by using the symplectic integrator
of type $\Sscr\Bscr\Ascr\Bscr_{3}\,$, which
is described in~\cite{Las-Rob-2001}. On the other hand, the secular
Hamiltonian model~\eqref{frm:hsec} has been integrated by using the
RK4 method; let us here recall that the comparisons would not have
changed significantly if we would have adopted the semi-analytic
computations described at the end of
subsection~\ref{sbs:semianalytical-results}. This is due to the fact
that the plots produced by applying the latter two methods superpose
in a practically indistinguishable way. Of course, the numerical
integrations have been started from two sets of initial conditions
that are coherent between them because both of them correspond to the
same values of the orbital elements, that are reported in
Table~\ref{tab:dopo-griglia}.  Figure~\ref{fig:secVScom} includes
those results that can be directly compared.  It clearly shows that
the secular model is able to substantially capture some of the main
dynamical features of the complete one: the quasi-periodicity of the
motions and the amplitudes of all the plotted quantities, that are the
eccentricities of the two planets and the difference of the arguments
of their pericenters. However, the comparisons that are reported in
Figure~\ref{fig:secVScom} are not satisfactory at all for what
concerns the fundamental angular velocities that generate all the
quasi-periodic motions. Indeed, the periods look rather
different. This is particularly true for what concerns the component
giving the slow modulation that can be easily observed in the orbital
laws induced by the flow of the complete planetary
Hamiltonian~\eqref{frm:ham3cpiniz}: none of the Fourier harmonics
composing the motion related to the secular model looks to play the
same role. We stress that the agreement is much worse with respect to
the one we get for the Sun--Jupiter--Saturn system. In that
case, the same kind of comparisons between the complete planetary
Hamiltonian and the secular one (still calculated at order two in the
masses) showed a minimal discrepancy, because the relative errors on
the angular velocities of the perihelion arguments of Jupiter and
Saturn are about 0.1\% and 1.3\%, respectively (see at the end of
Section~2 of~\cite{Loc-Gio-2000}).

\section{Conclusions and perspectives}
\label{sec:conclu}

Computer-assisted algorithms based on the construction of the
Kolmogorov normal form can provide very good performances when a toy
model is considered. For instance, in~\cite{Cel-Gio-Loc-2000} the
Authors rigorously prove the existence of invariant tori for values of
the perturbing terms that are relatively close to the expected
breakdown threshold. When applied to a model that is more similar to a
system of physical interest though, the same approach is able to prove
the KAM theorem only in situations that are still rather far from the
edge (see~\cite{Val-Loc-2021}).  Therefore, the results described in
the present paper can be considered astonishingly good, considering
the adherence of the model to the physical problem. On one hand, the
$\upsilon$~Andromed{\ae} planetary system is considered to be much
more perturbed with respect to the Sun--Jupiter--Saturn one studied,
for example, in~\cite{Loc-Gio-2000}. On the other hand, for both these
three-body problems, the computer-assisted proofs can produce the same
result\footnote{Let us recall that it is not yet available a
  computer-assisted proof of existence of invariant tori for the
  complete planetary Hamiltonian describing the three-body dynamics of
  the Sun--Jupiter--Saturn system.} when they are used at their best:
they ensure the existence of invariant tori into the framework of a
realistic secular model that is calculated at order two in the masses,
with the purpose to well approximate their planetary dynamics. A key
ingredient that we introduce in the present work is the construction
of a preliminary normal form that allows to approximate a lower
dimensional elliptic invariant torus. It should be noted that few
details introduced in the approach presented in this work lead to some
simplifications.  For instance, we have considered the initial
conditions reported in Table~\ref{tab:dopo-griglia} that are expected
to be the most robust among those compatible with the observational
data reported in~\cite{McArt-et-al-2010}. Moreover, the proof of
convergence relates only to the final normalization algorithm, instead
of the whole sequence of algorithms applied. It would be desirable to
further refine the techniques adopted in the context of this kind of
computer-assisted proofs, in order to apply them to realistic models
without any extra simplification, as it already happens for toy models
(see~\cite{Fig-Har-Luq-2017}).

In our view, the interest of our work goes beyond the rigorous
mathematical proofs ensuring the existence of invariant tori. The
adaptation of the constructive algorithm to the case of librational
tori emphasizes their role in stabilizing the planetary orbits. This
is strictly related to the criterion of robustness that allows us to
make selections among the initial conditions that are compatible with
the observations. Such a criterion will be explained in a forthcoming
publication and it will require a further extension to the systems in
Mean Motion Resonance (MMR). For such a reason, we stress the
importance of adapting our approach to MMR planetary models so as to
combine the constructions of the normal forms for both elliptic and
KAM tori. In the framework of MMR systems, the librational orbits are
more complicated with respect to those considered in the present
paper, but the effectiveness of the computational procedure
constructing an average model at order two in the masses has already
been shown (see~\cite{San-Lib-2019}).

The main advantage of normal form methods is that they allow to
describe the dynamics in a region of the phase space, instead of
analyzing one trajectory at a time, thus obtaining a result in a set of
zero measure. In particular, in this work we have designed an approach
based on a preliminary construction of the elliptic tori, before
studying the librational orbits in their neighborhood. Indeed, the
secular model here introduced shows some discrepancies between the
semi-analytical results and the motions generated by the complete
planetary Hamiltonian (as shown in Figure~\ref{fig:secVScom}). A
better agreement with the numerical solutions can be achieved by
extending the normal form approach from the secular Hamiltonian to the
complete one. In other words, we should apply the same adaption of the
constructive algorithm that was successfully implemented by passing
from~\cite{Loc-Gio-2000} to~\cite{Loc-Gio-2007} in the case of the
Sun--Jupiter--Saturn system.

The characterization of the three dimensional architecture of
extrasolar planetary systems is one of the very exciting problems that
is currently investigated in Dynamical Astronomy.  It has been studied
within the framework of both complete and secular planetary
Hamiltonians (see, e.g., \cite{Mich-FerM-Beau-2006},
\cite{Lib-Tsi-2009}
and~\cite{Vol-Roi-Lib-2019}). In~\cite{Vol-Loc-San-2018}, some of the
Authors of this article introduced a reverse approach based on KAM
theory: values of the mutual inclinations are considered to be
acceptable if the subsequent construction of the Kolmogorov normal
form is convergent. The effectiveness of that method was limited to
extrasolar systems hosting planets with very moderate eccentricities
(i.e., smaller than~$0.1$). The original motivation of the present
work was aiming to overcome this artificial constraint due to the
technique adopted in~\cite{Vol-Loc-San-2018}. In this respect, the
present work is certainly successful, because the average eccentricity
for stable planetary orbits is about~$0.3$ in the case of the model we
considered here for studying the $\upsilon$~Andromed{\ae} system.  In
a near future, we plan to restart the explorations of the~3D
architecture of extrasolar planetary systems, with the prescription of
the KAM stability of those models. Of course, this will be made by
exporting the new strategy introduced here, namely by combining the
constructions of the normal forms for both elliptic tori and KAM ones.

\section*{Data Availability}
The data that support the findings of this study are
available from the corresponding author, C.~C., upon
reasonable request.

\subsection*{Acknowledgments}
This work was partially supported by the MIUR-PRIN project 20178CJA2B
-- ``New Frontiers of Celestial Mechanics: theory and Applications'',
the ``Beyond Borders'' program of the University of Rome Tor Vergata
through the project ASTRID (CUP E84I19002250005), the ``Progetto
Giovani 2019'' program of the National Group of Mathematical Physics
(GNFM--INdAM) and the ASI Contract n. 2018-25-HH.0 (Scientific
Activities for JUICE, C/D phase).  The authors acknowledge also the
MIUR Excellence Department Project awarded to the Department of
Mathematics of the University of Rome ``Tor Vergata'' (CUP
E83C18000100006).

\appendix
\section{Computer-assisted
  proof of the existence of KAM tori: technicalities}
\label{sec:technicalities}
In this appendix we describe how to iteratively define the terms
making part of $\Sscr^{(r)}$, as it is defined
in~\eqref{frm:hamCAP}. Such an iterative scheme of estimates is
effectively implemented in a software package dealing with the
computer-assisted proof whose result is discussed in
section~\ref{sbs:rigorous-results}. Let us recall that this package is
publicly available from the web~\cite{Loc-CAP4KAM-2021} and it
subdivides the expansions in Hamiltonian terms according to slightly
different definitions with respect to those adopted in
subsection~\ref{sec:kamtor}. This is made in order to improve the
results obtained by performing the computer-assisted proof.  Indeed,
since the very beginning (i.e., when the normalization step $r=0$),
the codes in the already mentioned software
package~\cite{Loc-CAP4KAM-2021} consider expansions of
Hamiltonians of the following type:
\begin{equation}
  \label{frm:ham-step-r-kamtor}
  \Kscr^{(r)} = E^{(r)} +\tilde{ \vet \omega}\cdot \vet p + \sum_{s\ge
    0} \sum_{\ell \ge 2}f_\ell^{(r, s)}(\vet p, \vet q) + \sum_{s\ge
    r+1} \sum_{\ell = 0}^1 f_\ell^{(r, s)}(\vet p, \vet q)\ ,
\end{equation}
where $f_\ell^{(r, s)}\in\Pgot_{\ell\,,\,2s}$ $\forall\ \ell,\,s$.
A generic function $g$ belongs to the class of functions $\Pgot_{\ell\,,\,2s}$
if its Taylor--Fourier expansion is such that
\begin{equation}
g( \vet p, \vet q) =
\sum_{{j_1+\ldots+j_n=\ell}}\ \sum_{\max_{i=1,\ldots,n} |k_i|\le 2s} \,\vet p^{\vet j}
\big[c_{\vet j,\vet k}\cos({\vet k}\cdot{ \vet q})
  + d_{\vet j,\vet k}\sin({\vet k}\cdot{ \vet q})\big] \ ,
\label{frm:esempio-g-in-Pgot}
\end{equation}
where $n$ denotes the number of degrees of freedom and the
coefficients $c_{\vet j,\vet k}\,,\,d_{\vet j,\vet k}\in\reali$; in
the previous formula, we have also adopted the multi-index notation,
i.e., $\vet p^{\vet j}=p_1^{j_1}\cdot\ldots\cdot p_n^{j_n}$. This
means that the only difference with respect to to the expansions
described during the discussion of the formal algorithm in
subsection~\ref{sec:kamtor} is due to the fact that here we consider
the $\ell_\infty$-norm of the Fourier harmonic $\vet k$ in order to
make the separation with respect to different classes of function,
instead of the $\ell_1$-norm $|\vet{k}|$, that has also been called
the trigonometrical degree corresponding to $\vet{k}$. Let us recall
that the functional norm whose corresponding values for some
generating functions are reported in the plots of
Figures~\ref{fig:normegen-ell}--\ref{fig:stime} is defined in such
a way that
$$
\|g\| =
\sum_{{j_1+\ldots+j_n=\ell}}\ \sum_{\max_{i=1,\ldots,n} |k_i|\le 2s} \,
\left(\big|c_{\vet j,\vet k}\big| + \big| d_{\vet j,\vet k}\big|\right)
$$
(see the corresponding
footnote${\ref{nota:def-norma}}\atop{\phantom{1}}$ in
subsection~\ref{subsec:CAP}).  Being the transformations of
coordinates defined by the Lie series, the estimates of the functional
norm $\|\cdot \|$ that are needed can be summarized by the following
statement.
\begin{lemma}
\label{lem:stime}
Let $g$ be a generic function belonging to the class of functions
$\Pgot_{\ell\,,\,sK}\,$, for some non-negative integers $\ell$ and
$s$, with $K\in \interi^+$. Let us consider the Lie series
$\exp\lie{\chi_1^{(r)}}g$ and $\exp\lie{\chi_2^{(r)}}g$, where the
generating functions are such that $\chi_1^{(r)}= X^{(r)}+
\vet\xi^{(r)}\cdot \vet q$, being $X^{(r)}\in\Pgot_{0\,,\,rK}\,$,
$\vet\xi^{(r)}\in\reali^n$, and $\chi_2^{(r)}\in\Pgot_{1\,,\,rK}\,$.
Therefore, the following estimates hold true for the summands
appearing in the expansions of those Lie series:
\begin{align}
\label{frm:stime-chi1}
\left\|\frac{1}{j!}\Lie{\chi_1^{(r)}}^j g\right\| \le & \binom{\ell}{j}
\left( \max_{i}\left\{\left\| \frac{\partial X^{(r)}}{\partial q_{i}}\right\|
\right\}+ \max_{i}\left\{\left|\xi_{i}^{(r)}\right|\right\} \right)^j \|g\|\ ,
\\
\left\|\frac{1}{j!}\Lie{\chi_2^{(r)}}^j g \right\| \le & \frac{1}{j!}
\prod_{i=1}^{j} \left[ \ell \max_{i}\left\{\left\|
  \frac{\partial \chi_2^{(r)}}{\partial q_{i}}\right\|\right\}
  + K(s+(i-1)r)\max_{i}\left\{\left\|\frac{\partial\chi_{2}^{(r)}}{\partial p_{i}}
  \right\|\right\} \right] \|g\|\ .
\label{frm:stime-chi2}
\end{align}
\end{lemma}
Since the proof can be easily reconstructed by applying repeatedly the
definition of the Poisson brackets, it is omitted. Let us recall that
we have always adopted $K=2$ for all the computations described in
Sections~\ref{sec:toro} and~\ref{sec:results}.

The statement of lemma~\ref{lem:stime} makes evident that we need to
bound the derivatives of the generating functions, in such a way to
determine four majorants $\Gscr_{11}^{(r)}$, $\Gscr_{12}^{(r)}$, $\Gscr_{21}^{(r)}$ and $\Gscr_{22}^{(r)}$ satisfying the inequalities $\max_i
\|\partial X^{(r)}/\partial q_i\| \le \Gscr_{11}^{(r)}$, $\max_i |
\xi_i^{(r)}| \le \Gscr_{12}^{(r)}$, $\max_i \|\partial
\chi_2^{(r)}/\partial q_i\| \le \Gscr_{21}^{(r)}$ and $\max_i
\|\partial \chi_2^{(r)}/\partial p_i\| \le \Gscr_{22}^{(r)}$
at the generic normalization step $r$. When it is not possible
to explicitly calculate those derivatives, because the expansions
of $\chi_1^{(r)}$ and $\chi_2^{(r)}$ are unknown, it is convenient
to define their upper bounds as follows:
\begin{align}
\label{G1}
\Gscr_{11}^{(r)} & \le
2r\left( \min_{0<\max_i|k_i|\le 2r} |\vet k \cdot \vet \omega|\right)^{-1}
\Fscr_0^{(r-1,r)}\ ,
\quad
\Gscr_{12}^{(r)} \le  (m^{(r)})^{-1}\langle \Fscr_1^{(r-1,r)}\rangle\ ,
\\
\Gscr_{21}^{(r)} & \le
2r\left( \min_{0<\max_i|k_i|\le 2r} |\vet k \cdot \vet \omega|\right)^{-1}
\Fscr_1^{(r-1,r)}\ ,
\quad
\Gscr_{22}^{(r)} \le
\left( \min_{0<\max_i|k_i|\le 2r} |\vet k \cdot \vet \omega|\right)^{-1}
\Fscr_1^{(r-1,r)}\ ,
\label{G2}
\end{align}
being $\Fscr_\ell^{(r,s)}$ an estimate of the norm of the generic term
$f_\ell^{(r,s)}$ and $m^{(r)}$ is such that $|C^{(r)}\cdot v|\ge
m^{(r)}|v|\ \forall \ v\in \reali^n$, where $C^{(r)}$ is the matrix
defined by the equation $\frac 1 2 C^{(r)} \vet p \cdot \vet p =
\sum_{s=0}^{r} \langle f_2^{(r,s)}\rangle$.  The estimates in
lemma~\ref{lem:stime}, combined with the formul{\ae} for the
re-definitions of the Hamiltonian terms in~\eqref{frm:dopochi1}
and~\eqref{frm:dopochi2}, provide the following outcomes:
\begin{equation}
\vcenter{\openup1\jot\halign{
 \hbox {\hfil $\displaystyle {#}$}
&\hbox {$\displaystyle {#}$\hfil}
&\hbox {$\displaystyle {#}$\hfil}\cr
  \hat \Fscr_{\ell -j}^{(r, jr+s)} & \hookleftarrow
  \binom \ell j \left(\Gscr_{11}^{(r)}+\Gscr_{12}^{(r)} \right)^j
  \Fscr_{\ell}^{(s, s)} \quad
  &\forall \ \ell \ge 2, \ 0\le s < r, \ 1 \le j \le \ell\ ,
  \cr
  \hat \Fscr_{\ell -j}^{(r, jr+s)} & \hookleftarrow
  \binom \ell j \left(\Gscr_{11}^{(r)}+\Gscr_{12}^{(r)} \right)^j
  \Fscr_{\ell}^{(r-1, s)} \quad
  &\forall \ \ell \ge 1, \ s \ge r,
  \ 1 \le j \le \ell
  \cr
}}
\label{ridef:hatF}
\end{equation}
(after having initially set
$\hat\Fscr_{\ell}^{(r,s)}=\Fscr_{\ell}^{(r-1, s)}$ $\forall \ \ell \ge
0, \ s \ge 0$) and
\begin{equation}
\vcenter{\openup1\jot\halign{
 \hbox {\hfil $\displaystyle {#}$}
&\hbox {$\displaystyle {#}$\hfil}
&\hbox {$\displaystyle {#}$\hfil}\cr
  \Fscr_{\ell}^{(r, jr+s)} & \hookleftarrow \frac{1}{j!}\prod_{i=0}^{j-1}
  \left[\ell \Gscr_{21}^{(r)}+K((j-1-i)r+s)\Gscr_{22}^{(r)} \right]
  \hat \Fscr_{\ell}^{(s, s)} \quad
  &\forall \ \ell \ge 2, \ 0\le s \le r, \ j \ge 1\ ,
  \cr
  \Fscr_{\ell}^{(r, jr+s)} & \hookleftarrow  \frac{1}{j!}\prod_{i=0}^{j-1}
  \left[\ell \Gscr_{21}^{(r)}+K((j-1-i)r+s)\Gscr_{22}^{(r)} \right]
  \hat \Fscr_{\ell}^{(r, s)}
  \cr
  & \qquad\qquad\qquad \forall\ \ell \ge 0, \  s > r, j \ge 1 \ {\rm or}
  \ \ell =1, \ s = r,  \ j \ge 1
  \cr
}}
\label{ridef:F}
\end{equation}
(after having initially set
$\Fscr_{\ell}^{(r,s)}=\hat\Fscr_{\ell}^{(r, s)}$ $\forall \ \ell \ge
0, \ s \ge 0$). Of course, after having completed all the
re-definitions prescribed by formula~\eqref{ridef:hatF}, we have to
impose that $E^{(r)}=E^{(r-1)}+\hat\Fscr_{0}^{(r,r)}$ and
$\hat\Fscr_{0}^{(r, r)}=0$; analogously, at the end of the
re-definitions~\eqref{ridef:F}, we have to set $\hat\Fscr_{1}^{(r,
  r)}=0$. This is made in order to update the energy value that at the
end of the procedure will correspond to that of the invariant KAM
torus and to take into account of the homological
equations~\eqref{eq:kolgen1}, \eqref{eq:kolgen2}
and~\eqref{frm:trasl}. Furthermore, the estimate ensuring that the
non-degeneracy condition is satisfied has to be updated thanks to the
following redefinition:
\begin{equation}
\label{newrho}
m^{(r+1)} = m^{(r)} - 2 \Fscr_2^{(r, r)}.
\end{equation}

Let us emphasize that all the estimates~\eqref{G1}--\eqref{newrho} for
the quantities $\Gscr_{11}^{(s)}$, $\Gscr_{12}^{(s)}$,
$\Gscr_{21}^{(s)}$, $\Gscr_{22}^{(s)}$, $\hat{\Fscr}_{\ell}^{(r,s)}$,
$\Fscr_{\ell}^{(r,s)}$ and $m^{(s)}$ are used just when $R_{\rm I} < s
\le R_{\rm II}$, while the estimates are made directly from the
expansions of the (generating) functions that are computed explicitly
for all $0 \le s \le R_{\rm I}$. Of course, this is made in order to
take advantage as much as possible from the explicit computations.  We
stress that Figure~\ref{fig:stime} actually includes the plots of the
finite sequence of values $\Gscr_{22}^{(r)}$ (as a function of the
normalization step $r$), that are computed exactly as it has been
explained in the present appendix.

In order to complete the iterative definition of the set
$\Sscr^{(r)}$, we have to determine the three parameters that appear
in formula~\eqref{frm:hamCAP} and are ruling the decay of the infinite
tail of terms making part of the expansions of the
Hamiltonian~$\Hscr^{(r)}$; they are $\Escr_{r}\,$, $a_{r}$ and
$\zeta_{r}\,$. These final settings are omitted because they can be
easily found in subsection~3.2.3 of~\cite{Val-Loc-2021}.

\end{document}